\documentclass{sn-jnl}
\usepackage{graphicx}%
\usepackage{multirow}%
\usepackage{amsmath,amssymb,amsfonts}%
\usepackage{amsthm}%
\usepackage{mathrsfs}%
\usepackage[title]{appendix}%
\usepackage{xcolor}%
\usepackage{textcomp}%
\usepackage{manyfoot}%
\usepackage{booktabs}%
\usepackage{algorithm}%
\usepackage{algorithmicx}%
\usepackage{algpseudocode}%
\usepackage{listings}%
\usepackage{geometry}

\theoremstyle{thmstyleone}%
\newtheorem{theorem}{Theorem}[section]%  meant for continuous numbers
\newtheorem{lemma}{Lemma}[section] %
\theoremstyle{thmstyletwo}%

\theoremstyle{thmstylethree}%

\begin{document}

\title[Article Title]{Oriented diameter of the complete tripartite graph  (II)}

%%=============================================================%%
%% GivenName	-> \fnm{Joergen W.}
%% Particle	-> \spfx{van der} -> surname prefix
%% FamilyName	-> \sur{Ploeg}
%% Suffix	-> \sfx{IV}
%% \author*[1,2]{\fnm{Joergen W.} \spfx{van der} \sur{Ploeg} 
%%  \sfx{IV}}\email{iauthor@gmail.com}
%%=============================================================%%

\author*[1]{\fnm{Jing} \sur{Liu}}\email{2022201091@buct.edu.cn}

\author[2]{\fnm{Hui} \sur{Zhou}}\email{zhouh06@qq.com, zhouhlzu06@126.com}
\equalcont{These authors contributed equally to this work.}

\affil[1,2]{\orgdiv{College of Mathematics and Physics}}

\affil[]{ \orgname{Beijing University of Chemical Technology}, \city{Beijing}, \postcode{100029},  \country{China}}

%\affil[3]{\orgdiv{Department}, \orgname{Organization}, \orgaddress{\street{Street}, \city{City}, \postcode{610101}, \state{State}, \country{Country}}}

%%==================================%%
%% Sample for unstructured abstract %%
%%==================================%%

\abstract{\normalfont For a graph $G$, let $\mathbb{D}(G)$ denote the set of all strong orientations of $G$, and the oriented diameter of $G$ is $f(G)=\min \{\operatorname{diam}(D) \mid D \in \mathbb{D}(G)\}$, which is the minimum value of the diameters $\operatorname{diam}(D)$ where $D \in \mathbb{D}(G)$.  In this paper, we determine the oriented diameter of  complete tripartite graphs $K(3,3, q)$ and $K(3,4, q)$, these are  special cases that arise  in determining the oriented diameter of  $K(3, p, q)$.
}


%%\pacs[JEL Classification]{D8, H51}


\maketitle

\section{Introduction}\label{sec1}

Let $G$ be a finite undirected simple connected graph with vertex set $V(G)$ and edge set $E(G)$. The orientation of a graph $G$ is a directed graph obtained from $G$ by assigning each of its edges in $G$ a direction. An orientation $D$ of $G$ is strong if for any two vertices $u, v$ in $D$, there exists a directed path from $u$ to $v$. For any $u, v \in V(G)$, the distance $d_G(u, v)$ denotes the length of a shortest $(u, v)$ - path in $G$, which is the number of edges in a shortest path connecting $u$ and $v$ in $G$. The diamater of $G$ is defined as $\operatorname{diam}(G)=\max \left\{d_G(u, v) \mid u, v \in V(G)\right\}$. An edge $e \in E(G)$ is called a bridge of a graph $G$ if the subgraph obtained by deleting the edge $e$ of the graph $G$ is disconnected. A graph is called bridreless if it has no bridge.

Robbins' one-way street theorem \cite{r1} proves that a connected graph has a strong orientation if and only if it is bridgeless. Boesch and Tindell \cite{r2} proposed the notion $f(G)$ in order to extend  Robbins' theorem \cite{r1}.

Let $G$ be a connected and bridreless graph, and $\mathbb{D}(G)$ be the set of all strong orientations of $G$. Define the oriented diameter of $G$ to be
$$
f(G)=\min \{\operatorname{diam}(D) \mid D \in \mathbb{D}(G)\}.
$$
 For an arbitrary connected graph $G$,  the problem of evaluting oriented diameter $f(G)$ is very difficult. In reality, Chv$\acute{a}$tal and Thomassen \cite{r3} demonstrate that the problem of deciding whether a graph admits an orientation of diameter two is NP-hard. Next, we will present some results on the oriented diameter that have been obtained in the literature.

Given any positive integers $n, p_1, p_2, \dots, p_n$, let $K_n$ denote the complete graph of order $n$, and $K\left(p_1, p_2, \dots, p_n\right)$ denote the complete $n$-partite graph having $p_i$ vertices in the $i$-th partite set $V_i$ for each $i=1,2, \dots, n$. Thus $K_n$ is isomorphic to $K\left(p_1, p_2, \dots, p_n\right)$ where $p_1=p_2= \dots =p_n=1$. The oriented diameter of complete graph $K_n$ was obtained by Boesch and Tindell \cite{r2}:
$$
f\left(K_n\right)= \begin{cases}2, & \text { if } n \geq 3 \text { and } n \neq 4 ; \\ 3, & \text { if } n=4 .\end{cases}
$$
The oriented diameter of complete bipartite graph $K(p, q)$ was given by
\v{S}olt\'{e}s \cite{r4}:
$$
f(K(p, q))=\left\{\begin{array}{l}
	3, \text { if } 2 \leq p \leq q \leq\binom{ p}{\left\lfloor\frac{p}{2}\right\rfloor} ; \\
	4, \text { if } q>\binom{p}{\left\lfloor\frac{p}{2}\right\rfloor} ;
\end{array}\right.
$$
where $\lfloor x\rfloor$ denotes the greatest integer not exceeding  $x$.  Let $n \geq 3$, Plesn\'{i}k \cite{r5}, Gutin \cite{r6}, and Koh and Tay \cite{r7} independently obtained the oriented diameter of complete $n$-partite graph, the specific result is as follows:
$$
2 \leq f\left(K\left(p_1, p_2, \cdots, p_n\right)\right) \leq 3 .
$$
Let $p \geq 2$ and $n \geq 3$, Koh and Tan \cite{r8} also obtained:
$$
\begin{gathered}
	f(K(\overbrace{p, p, \cdots p}^n))=2 ;\\
	f(K(\overbrace{p, p, \cdots p}^r, q))=2, (r\geqslant 3, p\geqslant 3, 1\leqslant q\leqslant 2p). 
\end{gathered}
$$
They also got some other results on complete multipartite graphs.

Actually, a problem was proposed by Koh and Tan \cite{r7} : ``given a complete multipartite graph $G=K\left(p_1, p_2, \cdots, p_n\right)$, classify it according to $f(G)=2$ or $f(G)=3$."  We endeavor to classify complete multipartite graphs according to the oriented diameter 2 or 3.  Koh and Tan \cite{r9} obtained:
$$
f(K(2, p, q))=2 \text { for } 2 \leq p \leq q \leq\binom{ p}{\left\lfloor\frac{p}{2}\right\rfloor}.
$$
and the authors and Rao and Zhang proved \cite{r10}  that
$$
f(K(2, p, q))=3 \text { for } q>\binom{p}{\left\lfloor\frac{p}{2}\right\rfloor}.
$$
So far, the oriented diameter of complete tripartite graph $K(2, p, q)$ is completely determined.

In this paper, we devote ourselves to deal with the case $K(3, p, q)$ for $p \geq 3$, and we first discuss special cases $K(3,3, q)$ and $K(3,4, q)$, the  results are as follows:
$$
\begin{aligned}
	& f(K(3,3, q))= \begin{cases}2, & \text { if } q \leq 6 ; \\
		3, & \text { if } q>6; \end{cases} \\
	& f(K(3,4, q))= \begin{cases}2, & \text { if } q \leq 11 ; \\
		3, & \text { if } q>11 .\end{cases}
\end{aligned}
$$
In a long manuscript, we also have  determined the oriented diameter of $K(3,p,q)$ for $p\geqslant 5$. Hence the problem posed by Koh and Tay for complete tripartite graphs $K(3, p, q)$ is completely solved.

\section{Preliminaries}\label{sec2}

For a digraph $D$, we denote its vertex set by $V(D)$. Take any $u, v \in V(D)$, the  distance $\partial_D(u, v)$ denotes the number of directed arcs in a shortest directed path from $u$ to $v$ in $D$. The diameter of $D$ is defined as $\operatorname{diam}(D)=\max \left\{\partial_D(u, v) \mid u, v \in V(D)\right\}$. For $u, v \in V(D), U, V \subseteq V(D)$ and $U \cap V=\varnothing$, if the direction is from $u$ to $v$, we write `$u \rightarrow v$'; if $a \rightarrow b$ for each $a \in U$ and each $b \in V$, we write `$U \rightarrow V$'; if $U=\{u\}$, we write `$u \rightarrow V$' for `$U \rightarrow V$';  if $V=\{v\}$, we write `$U \rightarrow v$' for `$U \rightarrow V$'. In addition, the set $N_D^{+}(u)=\{x \in V(D) \mid u \rightarrow x\}$, which is a collection of all out-neighbors of $u$; the set $N_D^{-}(v)=\{y \in V(D) \mid y \rightarrow v\}$, which is a collection of all in-neighbors of $v$.

For the complete tripartite graph $K(3,p,q)$, $p\in \{3,4\}$, let
$$
\begin{aligned}
	& V_1=\left\{x_1, x_2, x_3\right\}, \\
	& V_2=\left\{y_1, y_2, \cdots, y_p\right\},  \\
	& V_3=\left\{z_1, z_2, \cdots, z_q\right\}.
\end{aligned}
$$
be the three parts of the vertex set of $K(3,p, q)$. Let $D$ be a strong orientation of $K(3,p,q)$. We consider the sets 
$$
\begin{aligned}
	& N_D^{+++}=N_D^{+}\left(x_1\right) \cap N_D^{+}\left(x_2\right) \cap N_D^{+}\left(x_3\right), \\
	& N_D^{++-}=N_D^{+}\left(x_1\right) \cap N_D^{+}\left(x_2\right) \cap N_D^{-}\left(x_3\right), \\
	& N_D^{+-+}=N_D^{+}\left(x_1\right) \cap N_D^{-}\left(x_2\right) \cap N_D^{+}\left(x_3\right), \\
	& N_D^{-++}=N_D^{-}\left(x_1\right) \cap N_D^{+}\left(x_2\right) \cap N_D^{+}\left(x_3\right), \\
	& N_D^{+--}=N_D^{+}\left(x_1\right) \cap N_D^{-}\left(x_2\right) \cap N_D^{-}\left(x_3\right), \\
	& N_D^{-+-}=N_D^{-}\left(x_1\right) \cap N_D^{+}\left(x_2\right) \cap N_D^{-}\left(x_3\right), \\
	& N_D^{--+}=N_D^{-}\left(x_1\right) \cap N_D^{-}\left(x_2\right) \cap N_D^{+}\left(x_3\right), \\
	& N_D^{---}=N_D^{-}\left(x_1\right) \cap N_D^{-}\left(x_2\right) \cap N_D^{-}\left(x_3\right),
\end{aligned}
$$
For $i \in\{2,3\}$, the following eight sets
$$
\begin{aligned}
	& V_i^{+++}=V_i \cap N_D^{+++}, \\
	& V_i^{++-}=V_i \cap N_D^{++-}, \\
	& V_i^{+-+}=V_i \cap N_D^{+-+}, \\
	& V_i^{-++}=V_i \cap N_D^{-++}, \\
	& V_i^{+--}=V_i \cap N_D^{+--}, \\
	& V_i^{-+-}=V_i \cap N_D^{-+-}, \\
	& V_i^{--+}=V_i \cap N_D^{--+}, \\
	& V_i^{---}=V_i \cap N_D^{---} .
\end{aligned}
$$
form a partition of $V_i$.  For convenience, we will denote $V_i^{+++}$as $V_i^{+}  and \   \ V_i^{---}$as $V_i^{-}$.\\

\begin{lemma}
	
	Suppose $\{i, j\}=\{2,3\}$.\\
	1. If $V_i^{+++} \neq \emptyset$, then $V_i^{+++} \rightarrow V_j$ and $\left|V_i^{+++}\right|=1$; if $V_i^{---} \neq \emptyset$, then $V_j \rightarrow V_i^{---}$and $\left|V_i^{---}\right|=1$.\\
	2. If $V_i^{+++} \neq \emptyset$, then $V_j^{+++}=\emptyset$; if $V_i^{---} \neq \emptyset$, then $V_j^{---}=\emptyset$.
	
\end{lemma}

\textit{\textbf{Proof.}} 
Suppose $V_i^{+++} \neq \emptyset$. Take any $y \in V_i^{+++}$and any $z \in V_j$, we have $\partial_D(y, z) \leqslant 2$. If $z \rightarrow y$, then $N_D^{+}(y) \subseteq V_j \backslash\{z\}$. We know $\partial_D\left(z^{\prime}, z\right) \geqslant 2$ for any $z^{\prime} \in V_j \backslash\{z\}$, so $\partial_D(y, z) \geqslant 3$, a contradiction. Hence $y \rightarrow z$. This means $V_i^{+++} \rightarrow V_j$.

   Suppose $X$ is a strongly connected digraph. Let $u, v \in V(X)$ be two vertices of $X$. If $N_X^{+}(u) \cap N_X^{-}(v)=\emptyset$, then $\partial_X(u, v) \neq 2$.  We assume $\partial_X(u, v)=2$, then there exists $w \in V(X)$ such that $u \rightarrow w \rightarrow v$. So $w \in N_X^{+}(u) \cap N_X^{-}(v) \neq \emptyset$, a contradiction.  For distinct vertices $y_h, y_k \in V_i^{+++}$, we have $N_D^{+}\left(y_h\right) \subseteq V_j$ and $N_D^{-}\left(y_k\right) \subseteq$ $V_1$. $\quad V_1 \cap V_j=\emptyset$ implies $N_D^{+}\left(y_h\right) \cap N_D^{-}\left(y_k\right)=\emptyset$, so we get $\partial_D\left(y_h, y_k\right) \geqslant 3$, a contradiction. Thus $\left|V_i^{+++}\right|=1$. The proof for the case $V_i^{---} \neq \emptyset$ is analogous.

Suppose $V_i^{+++} \neq \emptyset$ and $V_j^{+++} \neq \emptyset$, then $V_i^{+++} \rightarrow V_j$ and $V_j^{+++} \rightarrow V_i$, i.e., for $y \in V_i^{+++}$and $z \in V_j^{+++}$, we have $y \rightarrow z$ and $z \rightarrow y$, a contradiction. The proof for the case $V_i^{---} \neq \emptyset$ is analogous.  $\hfill\blacksquare$ \\

\section{ The oriented diameter of $K(3,3,q)$ }\label{sec3}

\begin{lemma}
	For $3 \leq q \leq 6,  \  f(K(3,3, q))=2$.

\end{lemma} 
 \textit{\textbf{Proof.}} 
When $q=3$, it follows from Theorem 3 in Koh and Tan \cite{r8} (Discrete
Mathematics, 1996) that $f(K(3,3, q))=2$.

When $q=6$, let $V_3=V_3^{+++} \cup V_3^{+--} \cup V_3^{-+-} \cup V_3^{---}$where $\left|V_3^{+++}\right|=\left|V_3^{---}\right|=1$ and $\left|V_3^{+--}\right|=\left|V_3^{-+-}\right|=2$. Orient $K(3,3,6)$ as follows : $\left\{y_2, y_3\right\} \rightarrow\left\{x_1, x_2\right\} \rightarrow y_1 \rightarrow x_3 \rightarrow\left\{y_2, y_3\right\}, V_3^{+++} \rightarrow y_1 \rightarrow V_3 \backslash V_3^{+++}$, $V_3^{+++} \rightarrow\left\{y_2, y_3\right\} \rightarrow V_3^{---}$and it is enough for $\left\{y_2, y_3\right\}$ to form a directed four-length circle with two vertices in $V_3^{+--}$or $V_3^{-+-}$. Let $D_{6}$ be the resulting orientation, it is easy to verify that $\operatorname{diam}(D_{6})=2$. On the basis of the above orientation $D_{6}$, if $V_3^{---}=\varnothing$,  then we can get $f(K(3,3,5))=2$.

When $q=4$, let $V_3=V_3^{+++} \cup V_3^{++-} \cup V_3^{+-+} \cup V_3^{+--}$where $\left|V_3^{+++}\right|=\left|V_3^{++-}\right|=\left|V_3^{+-+}\right|=\left|V_3^{+--}\right|=1$. Orient $K(3,3,4)$ as follows: $V_2 \rightarrow x_1$, $\left\{y_2, y_3\right\} \rightarrow x_2 \rightarrow y_1,\left\{y_1, y_3\right\} \rightarrow x_3 \rightarrow y_2, V_3^{+++} \cup V_3^{++-} \rightarrow y_1 \rightarrow V_3^{+-+} \cup V_3^{+--}$, $V_3^{+++} \cup V_3^{+-+} \rightarrow y_2 \rightarrow V_3^{++-} \cup V_3^{+--}$and $V_3 \rightarrow y_3$. Let $D_{4}$ be the resulting orientation, it is easy to verify that $\operatorname{diam}(D_{4})=2$.  $\hfill\blacksquare$ \\

\begin{theorem}
	Suppose $q \geq 3$, then
	$$
	f(K(3,3, q))= \begin{cases}2, & \text { if } q \leq 6 ; \\ 3, & \text { if } q>6.  \end{cases}
	$$
	%\begin{eqnarray}\label{10}
	%....
	%\end{eqnarray}
\end{theorem}

\textit{\textbf{Proof.}} When $3 \leq q \leq 6$, we have shown in Lemma 3.1 that there exists an orientation of diameter 2 for $K(3,3, q)$.  When $q>6$, we prove it by contradiction. Assuming $f(K(3,3, q))=2$ when $q>6$, then $K(3,3, q)$ has a strong orientation $D$ with diameter $\operatorname{diam}(D)=2$.  Let $V_1=\left\{x_1, x_2, x_3\right\}, V_2=\left\{y_1, y_2, y_3\right\}$ and $V_3=\left\{z_1, z_2, \cdots ,z_q\right\}$ be the three parts of the vertex set of $K(3,3, q)$, and $i=\left|N_D^{+}\left(x_1\right) \cap V_2\right|, j=\left|N_D^{+}\left(x_2\right) \cap V_2\right|, k=\left|N_D^{+}\left(x_3\right) \cap V_2\right|$. For any $i, j, k \in\{0,1,2,3\}$, we have the following ten exhaustive cases Case(i,j,k)s to obtain the required contradiction.\\
(1) Case $(0,0,0)$. $ V_2 \rightarrow V_1$.

Take any $x \in V_1, z \in V_3$, since $\partial_D(x, z) \leq 2$, then we have $x \rightarrow z$, this means $V_1 \rightarrow V_3$. So $\partial_D\left(x_1, x_2\right) \geq 3$, a contradiction.\\
(2) Case $(0,0,1)$. $V_2 \rightarrow\left\{x_1, x_2\right\},\left\{y_2, y_3\right\} \rightarrow x_3 \rightarrow y_1$.

Take any $x \in V_1 \backslash\left\{x_3\right\}$ and $z \in V_3$, since $\partial_D(x, z) \leq 2$, then we have $x \rightarrow z$, this means $V_1 \backslash\left\{x_3\right\} \rightarrow V_3$. So $\partial_D\left(x_1, x_2\right) \geq 3$, a contradiction.\\
(3) Case $(0,0,2)$. $V_2 \rightarrow\left\{x_1, x_2\right\}, y_3 \rightarrow x_3 \rightarrow\left\{y_1, y_2\right\}$.

Take any $x \in V_1 \backslash\left\{x_3\right\}$ and $z \in V_3$, since $\partial_D(x, z) \leq 2$, then we have $x \rightarrow z$, this means $V_1 \backslash\left\{x_3\right\} \rightarrow V_3$. So $\partial_D\left(x_1, x_2\right) \geq 3$, a contradiction.\\
(4) Case $(0,0,3)$. $V_2 \rightarrow\left\{x_1, x_2\right\}, x_3 \rightarrow V_2$.

Take any $x \in V_1 \backslash\left\{x_3\right\}$ and $z \in V_3$, since $\partial_D(x, z) \leq 2$, then we have $x \rightarrow z$, this means $V_1 \backslash\left\{x_3\right\} \rightarrow V_3$. So $\partial_D\left(x_1, x_2\right) \geq 3$, a contradiction.\\
(5) Case $(0,1,1)$.\\
Subcase 1: $V_2 \rightarrow x_1,\left\{y_2, y_3\right\} \rightarrow\left\{x_2, x_3\right\} \rightarrow y_1$.

Take any $y \in V_2 \backslash\left\{y_1\right\}$ and $z \in V_3$, since $\partial_D(z, y) \leq 2$, then we have $z \rightarrow y$, this means $V_3 \rightarrow V_2 \backslash\left\{y_1\right\}$. So $\partial_D\left(y_2, y_3\right) \geq 3$, a contradiction.\\
Subcase2: $V_2 \rightarrow x_1,\left\{y_2, y_3\right\} \rightarrow x_2 \rightarrow y_1,\left\{y_1, y_3\right\} \rightarrow x_3 \rightarrow y_2$.

We know $y_1 \in V_2^{-+-}, y_2 \in V_2^{--+}$and $y_3 \in V_2^{---}$. By Lemma 2.1, we can get $V_3^{-+-}=V_3^{--+}=V_3^{---}=\varnothing$. Take any $z \in V_3$, since $\partial_D\left(x_1, z\right) \leq 2$ and $\partial_D\left(z, y_3\right) \leq 2$, then we have $x_1 \rightarrow z$ and $z \rightarrow y_3$, this means $x_1 \rightarrow V_3$ and $V_3 \rightarrow y_3$. Thus $V_3=V_3^{+++} \cup V_3^{++-} \cup V_3^{+-+} \cup V_3^{+--}$. Since $\partial_D\left(V_3^{++-}, y_i\right) \leq 2$ \ where $i=1,2,3$, \ then we have $V_3^{++-} \rightarrow y_1 . \ $ Since $\partial_D\left(V_3^{+-+}, y_i\right) \leq 2$\ where $i=1,2,3$,\ then we have $V_3^{+-+} \rightarrow y_2$. \ Since $\partial_D\left(x_3, V_3^{++-}\right) \leq 2, \partial_D\left(x_2, V_3^{+-+}\right) \leq 2, \partial_D\left(x_2, V_3^{+--}\right) \leq 2$ and $\partial_D\left(x_3, V_3^{+--}\right) \leq 2$, then we have $y_2 \rightarrow V_3^{++-}, y_1 \rightarrow V_3^{+-+}$and $\left\{y_1, y_2\right\} \rightarrow V_3^{+--}$. \ Let $\left|V_3^{++-}\right|=q_1,\left|V_3^{+-+}\right|=q_2$ and $\left|V_3^{+--}\right|=q_3$. \ If $q_1 \geq 2$, then there exists $z_i, z_j \in V_3^{++-}$such that $\partial_D\left(z_i, z_j\right) \geq 3$, a contradiction. \ Hence $q_1 \leq 1$. \ If $q_2 \geq 2$, then there exists $z_i, z_j \in V_3^{+-+}$such that $\partial_D\left(z_i, z_j\right) \geq 3$, a contradiction. \ Hence $q_2 \leq 1$. If $q_3 \geq 2$, then there exists $z_i, z_j \in V_3^{+--}$such that $\partial_D\left(z_i, z_j\right) \geq 3$, a contradiction. \ Hence $q_3 \leq 1$. \ By Lemma 2.1, we also have $\left|V_3^{+++}\right| \leq 1$, \ thus $\left|V_3\right|=q \leq 1+1+1+1=4 \leq 6$.\\
(6) Case $(0,1,2)$\\
Subcase 1: $V_2 \rightarrow x_1,\left\{y_2, y_3\right\} \rightarrow x_2 \rightarrow y_1, y_3 \rightarrow x_3 \rightarrow\left\{y_1, y_2\right\}$.

We know $y_1 \in V_2^{-++}, y_2 \in V_2^{--+}$and $y_3 \in V_2^{---}$. By Lemma 2.1, we can get $V_3^{-++}=V_3^{--+}=V_3^{---}=\varnothing$. Take any $z \in V_3$, since $\partial_D\left(x_1, z\right) \leq 2$ and $\partial_D\left(z, y_3\right) \leq 2$, then we have $x_1 \rightarrow z$ and $z \rightarrow y_3$, this means $x_1 \rightarrow V_3$ and $V_3 \rightarrow y_3$. Thus $V_3=V_3^{+++} \cup V_3^{++-} \cup V_3^{+-+} \cup V_3^{+--}$. Since $\partial_D\left(V_3^{+-+}, y_i\right) \leq 2$ where $i=1,2,3$, then we have $V_3^{+-+} \rightarrow y_2$. Since $\partial_D\left(x_2, V_3^{+-+}\right) \leq 2$ and $\partial_D\left(x_2, V_3^{+--}\right) \leq 2$, we have $y_1 \rightarrow V_3^{+-+}$ and $y_1 \rightarrow V_3^{+--}$. Let $\left|V_3^{++-}\right|=q_1,\left|V_3^{+-+}\right|=q_2,\left|V_3^{+--}\right|=q_3$ and $F=D\left[V_2 \backslash\left\{y_3\right\} \cup V_3^{++-}\right]$. Then $F$ is an orientation of $K\left(2, q_1\right)$ where $q_1 \leq q$. If $q_1 \geq 3$, then there exists $z_i, z_j \in V_3^{++-}$such that $\partial_F\left(z_i, z_j\right)=4$,\ so $\partial_D\left(z_i, z_j\right) \geq 3$, a contradiction. Hence $q_1 \leq 2$. If $q_2 \geq 2$, then there exists $z_i, z_j \in V_3^{+-+}$such that $\partial_D\left(z_i, z_j\right) \geq 3$, a contradiction. Hence $q_2 \leq 1$. If $q_3 \geq 2$, then there exists $z_i, z_j \in V_3^{+--}$such that $\partial_D\left(z_i, z_j\right) \geq 3$, a contradiction. Hence $q_3 \leq 1$. By Lemma 2.1,\ we also have $\left|V_3^{+++}\right| \leq 1$, thus $\left|V_3\right|=q \leq 1+2+1+1=5 \leq 6$.\\
Subcase 2: $V_2 \rightarrow x_1,\left\{y_2, y_3\right\} \rightarrow x_2 \rightarrow y_1 \rightarrow x_3 \rightarrow\left\{y_2, y_3\right\}$.

We know $y_1 \in V_2^{-+-}$and $y_2, y_3 \in V_2^{--+}$. By Lemma 2.1, we can get $V_3^{--+}=V_3^{-+-}=\varnothing$. Take any $z \in V_3$, since $\partial_D\left(x_1, z\right) \leq 2$, then we have $x_1 \rightarrow z$, this means $x_1 \rightarrow V_3$. Thus $V_3=V_3^{+++} \cup V_3^{++-} \cup V_3^{+-+} \cup V_3^{+--}$. Since $\partial_D\left(V_3^{++-}, y_i\right) \leq 2$ where $i=1,2,3$, then we have $V_3^{++-} \rightarrow y_1$. Since $\partial_D\left(V_3^{+-+}, y_i\right) \leq 2$ where $i=1,2,3$, then we have $V_3^{+-+} \rightarrow\left\{y_2, y_3\right\} . \quad$ Since $\partial_D\left(x_2, V_3^{+-+}\right) \leq 2$ and $\partial_D\left(x_2, V_3^{+--}\right) \leq 2$, then we have $y_1 \rightarrow V_3^{+-+}$and $y_1 \rightarrow V_3^{+--}$. \ If $V_3^{+-+} \neq \varnothing$, \ then $\partial_D\left(V_3^{+-+}, x_3\right) \geq 3$, a contradiction. So $V_3^{+-+}=\varnothing$. Let $\left|V_3^{++-}\right|=q_1,\left|V_3^{+--}\right|=q_2$ and $F_1=D\left[V_2^{--+} \cup V_3^{++-}\right], F_2=D\left[V_2^{--+} \cup V_3^{+--}\right]$. Then $F_1$ and $F_2$ are respextively an orientation of $K\left(2, q_1\right)$ and $K\left(2, q_2\right)$ where $q_1, q_2 \leq q$. \ If $q_1 \geq 3$,\ then there exists $z_i, z_j \in V_3^{++-}$such that $\partial_{F_1}\left(z_i, z_j\right)=4$, so $\partial_D\left(z_i, z_j\right) \geq 3$,\ a contradiction. \ Hence $q_1 \leq 2$. \ If $q_2 \geq 3$,\ then there exists $z_i, z_j \in V_3^{+--}$such that $\partial_{F_2}\left(z_i, z_j\right)=4$,\ so $\partial_D\left(z_i, z_j\right) \geq 3,\ $ a contradiction. \ Hence $q_2 \leq 2$. \ By Lemma 2.1,\ we also have $\left|V_3^{+++}\right| \leq 1$,\ thus $\left|V_3\right|=q \leq 1+2+2=5 \leq 6$.\\
(7) Case(0,1,3). $V_2 \rightarrow x_1,\left\{y_2, y_3\right\} \rightarrow x_2 \rightarrow y_1, x_3 \rightarrow V_2$.

We know $y_1 \in V_2^{-++}$and $y_2, y_3 \in V_2^{--+}$. \ By Lemma 2.1, we can get $V_3^{-++}=V_3^{--+}=\varnothing$. Take any $z \in V_3$, since $\partial_D\left(x_1, z\right) \leq 2$ and $\partial_D\left(z, x_3\right) \leq 2$,\ then we have $x_1 \rightarrow z$ and $z \rightarrow x_3$,\ this means $x_1 \rightarrow V_3$ and $V_3 \rightarrow x_3$. \ Thus $V_3=V_3^{++-} \cup V_3^{+--}$. Since $\partial_D\left(x_2, V_3^{+--}\right) \leq 2$,\ then we have $y_1 \rightarrow V_3^{+--}$. Let $\left|V_3^{++-}\right|=q_1\left|V_3^{+--}\right|=q_2$ and $F_1=\left[V_2 \cup V_3^{++-}\right], F_2=\left[V_2^{--+} \cup V_3^{+--}\right]$. Then $F_1$ and $F_2$ are respectively an orientation of $K\left(2, q_1\right)$ and $K\left(2, q_2\right)$ where $q_1, q_2 \leq q$. \ If $q_1 \geq 4$, then there exists $z_i, z_j \in V_3^{++-}$such that $\partial_{F_1}\left(z_i, z_j\right)=4$, so $\partial_D\left(z_i, z_j\right) \geq 3$, a contradiction. Hence $q_1 \leq 3$. \ If $q_2 \geq 3$, then there exists $z_i, z_j \in V_3^{+--}$such that $\partial_{F_2}\left(z_i, z_j\right)=4$, so $\partial_D\left(z_i, z_j\right) \geq 3$, a contradiction. Hence $q_2 \leq 2$. \ Thus $\left|V_3\right|=q \leq 3+2=5 \leq 6$.\\
(8) Case $(0,2,2)$\\
Subcase 1: $V_2 \rightarrow x_1, y_3 \rightarrow\left\{x_2, x_3\right\} \rightarrow\left\{y_1, y_2\right\}$.

We know $y_1, y_2 \in V_2^{-++}$and $y_3 \in V_2^{---}$. \ By Lemma 2.1, we can get $V_3^{-++}=V_3^{---}=\varnothing$. Take any $z \in V_3$, since $\partial_D\left(x_1, z\right) \leq 2$ and $\partial_D\left(z, y_3\right) \leq 2$,\ then we have $x_1 \rightarrow z$ and $z \rightarrow y_3$, this means $x_1 \rightarrow V_3$ and $V_3 \rightarrow y_3$. \ Thus $V_3=V_3^{+++} \cup V_3^{++-} \cup V_3^{+-+} \cup V_3^{+--}$. \ Let $\left|V_3^{++-}\right|=q_1,\left|V_3^{+-+}\right|=q_2,\left|V_3^{+--}\right|=q_3$ and $F_1=D\left[V_2^{-++} \cup V_3^{++-}\right], F_2=D\left[V_2^{-++} \cup V_3^{+-+}\right], F_3=D\left[V_2^{-++} \cup V_3^{+--}\right]$. Then $F_1, F_2$ and $F_3$ are respectively an orientation of $K\left(2, q_1\right), K\left(2, q_2\right)$ and $K\left(2, q_3\right)$ where $q_1, q_2, q_3 \leq q$. \ If $q_1 \geq 3$, then there exists $z_i, z_j \in V_3^{++-}$such that $\partial_{F_1}\left(z_i, z_j\right)=4$, so $\partial_D\left(z_i, z_j\right) \geq 3$, a contradiction. Hence $q_1 \leq 2$. \ If $q_2 \geq 3$, then there exists $z_i, z_j \in V_3^{+-+}$such that $\partial_{F_2}\left(z_i, z_j\right)=4$, so $\partial_D\left(z_i, z_j\right) \geq 3$, a contradiction. Hence $q_2 \leq 2$. \ If $q_3 \geq 3$, then there exists $z_i, z_j \in V_3^{+--}$such that $\partial_{F_3}\left(z_i, z_j\right)=4$, so $\partial_D\left(z_i, z_j\right) \geq 3$, a contradiction. Hence $q_3 \leq 2$. Since $\partial_D\left(V_3^{++-}, V_3^{+--}\right) \leq 2$, \ if $q_1=q_3=2$, then there exists $z_i \in V_3^{++-}, z_j \in V_3^{+--}$such that $\partial_D\left(z_i, z_j\right) \geq 3$, a contradiction. Hence we have $q_1 \leq 1$ or $q_3 \leq 1$. \ The argument for these two cases are similar, so we may assume $q_1 \leq 1$. \ By Lemma 2.1,\ we also have $\left|V_3^{+++}\right| \leq 1$, \ thus $\left|V_3\right|=q \leq 1+1+2+2=6 \leq 6$.\\
Subcase2: $V_2 \rightarrow x_1, y_3 \rightarrow x_2 \rightarrow\left\{y_1, y_2\right\}, y_1 \rightarrow x_3 \rightarrow\left\{y_2, y_3\right\}$.

We know $y_1 \in V_2^{-+-}, y_2 \in V_2^{-++}$and $y_3 \in V_2^{--+}$. \ By Lemma 2.1, we can get $V_3^{-+-}=V_3^{-++}=V_3^{--+}=\varnothing$. \ Take any $z \in V_3$, since $\partial_D\left(x_1, z\right) \leq 2$, then we have  $x_1 \rightarrow z$, this means $x_1 \rightarrow V_3$. \ Thus $V_3=V_3^{+++} \cup V_3^{++-} \cup V_3^{+-+} \cup V_3^{+--}$. Since $\partial_D\left(V_3^{++-}, y_i\right) \leq 2$  where $i=1,2,3$, then we have $V_3^{++-} \rightarrow y_1$. Since $\partial_D\left(V_3^{+-+}, y_i\right) \leq 2$ where $i=1,2,3$, then we have $V_3^{+-+} \rightarrow y_3$. Since $\partial_D\left(V_3^{++-}, x_2\right) \leq 2, \partial_D\left(x_3, V_3^{++-}\right) \leq 2, \partial_D\left(V_3^{+-+}, x_3\right) \leq 2$ and $\partial_D\left(x_2, V_3^{+-+}\right) \leq 2$, then we have $y_2 \rightarrow V_3^{++-} \rightarrow y_3$ and $y_2 \rightarrow V_3^{+-+} \rightarrow y_1$. Let $\left|V_3^{++-}\right|=q_1,\left|V_3^{+-+}\right|=q_2,\left|V_3^{+---}\right|=q_3$ and $F=D\left[V_2 \cup V_3^{+--}\right]$. Then $F$ is an orientation of $K\left(2, q_3\right)$ where $q_3 \leq q$. \ If $q_3 \geq 4$, then there exists $z_i, z_j \in V_3^{+--}$such that $\partial_F\left(z_i, z_j\right)=4$, so $\partial_D\left(z_i, z_j\right) \geq 3$, a contradiction.  \ Hence $q_3 \leq 3$. \ If $q_1 \geq 2$, then there exists $z_i, z_j \in V_3^{++-}$such that $\partial_D\left(z_i, z_j\right) \geq 3$, a contradiction. \ Hence $q_1 \leq 1$. \ If $q_2 \geq 2$, then there exists $z_i, z_j \in V_3^{+-+}$such that $\partial_D\left(z_i, z_j\right) \geq 3$, a contradiction. \ Hence $q_2 \leq 1$. \ By Lemma 2.1,\ we also have $\left|V_3^{+++}\right| \leq 1$, \ thus $\left|V_3\right|=q \leq 1+1+1+3=6 \leq 6$.\\
(9) Case $(1,1,1)$.\\
Subcase1: $V_1 \rightarrow y_1,\left\{y_2, y_3\right\} \rightarrow V_1$.

Take any $y \in V_2 \backslash\left\{y_1\right\}$ and $z \in V_3$, since $\partial_D(z, y) \leq 2$,\ then we have $z \rightarrow y$,\ this means $V_3 \rightarrow\left\{y_2, y_3\right\}$. \ So $\partial_D\left(y_2, y_3\right) \geq 3$,\ a contradiction.\\
Subcase2: $\left\{y_2, y_3\right\} \rightarrow\left\{x_1, x_2\right\} \rightarrow y_1,\left\{y_1, y_3\right\} \rightarrow x_3 \rightarrow y_2$.

We know $y_1 \in V_2^{++-}, y_2 \in V_2^{--+}$and $y_3 \in V_2^{---}$. \ By Lemma 2.1, we can get $V_3^{++-}=V_3^{--+}=V_3^{---}=\varnothing$. \ Thus $V_3=V_3^{+++} \cup V_3^{+-+} \cup V_3^{-++} \cup V_3^{+--} \cup V_3^{-+-}$. \ Take any $z \in V_3$,\ since $\partial_D\left(z, y_3\right) \leq 2$,\ then we have $z \rightarrow y_3$, this means $V_3 \rightarrow y_3$. Since $\partial_D\left(V_3^{+-+}, y_i\right) \leq 2$ where $i=1,2,3$, then we have $V_3^{+-+} \rightarrow y_2$. Since $\partial_D\left(V_3^{-++}, y_i\right) \leq 2$ where $i=1,2,3$, then we have $V_3^{-++} \rightarrow y_2$. Since $\partial_D\left(y_i, V_3^{+--}\right) \leq 2$ where $i=1,2,3$, then we have $y_1 \rightarrow V_3^{+--}$. Since $\partial_D\left(y_i, V_3^{-+-}\right) \leq 2$ where $i=1,2,3$, then we have $y_1 \rightarrow V_3^{-+-}$. Since $\partial_D\left(x_2, V_3^{+-+}\right) \leq 2, \partial_D\left(x_1, V_3^{-++}\right) \leq 2, \partial_D\left(x_3, V_3^{+--}\right) \leq 2$ and $\partial_D\left(x_3, V_3^{-+-}\right) \leq 2$, then we have $y_1 \rightarrow V_3^{+-+}, y_1 \rightarrow V_3^{-++}, y_2 \rightarrow V_3^{+--}$and $y_2 \rightarrow V_3^{-+-}$. \ Let $\left|V_3^{+-+}\right|=q_1,\left|V_3^{-++}\right|=q_2,\left|V_3^{+--}\right|=q_3$ and $\left|V_3^{-+-}\right|=q_4$. \ If $q_1 \geq 2$, then there exists $z_i, z_j \in V_3^{++-}$such that $\partial_D\left(z_i, z_j\right) \geq 3$, a contradiction. Hence $q_1 \leq 1$. \ If $q_2 \geq 2$, then there exists $z_i, z_j \in V_3^{-++}$such that $\partial_D\left(z_i, z_j\right) \geq 3$, a contradiction. Hence $q_2 \leq 1$. \ If $q_3 \geq 2$, then there exists $z_i, z_j \in V_3^{+--}$such that $\partial_D\left(z_i, z_j\right) \geq 3$, a contradiction. Hence $q_3 \leq 1$.  \ If $q_4 \geq 2$, then there exists $z_i, z_j \in V_3^{-+-}$such that $\partial_D\left(z_i, z_j\right) \geq 3$, a contradiction. Hence $q_1 \leq 1$. \ By Lemma 2.1, we also have $\left|V_3^{+++}\right| \leq 1$, thus $\left|V_3\right|=q \leq 1+1+1+1+1=5 \leq 6$.\\
Subcase3: $\left\{y_2, y_3\right\} \rightarrow x_1 \rightarrow y_1,\left\{y_1, y_3\right\} \rightarrow x_2 \rightarrow y_2,\left\{y_1, y_2\right\} \rightarrow x_3 \rightarrow y_3$.

We know $y_1 \in V_2^{+--}, y_2 \in V_2^{-+-}$and $y_3 \in V_2^{--+}$. \ By Lemma 2.1, we can get
$V_3^{+--}=V_3^{-+-}=V_3^{--+}=\varnothing$, thus $V_3=V_3^{+++} \cup V_3^{++-} \cup V_3^{+-+} \cup V_3^{-++} \cup V_3^{---}$. \ Since $\partial_D\left(V_3^{++-}, y_i\right) \leq 2$ where $i=1,2,3$, then we have $V_3^{++-} \rightarrow\left\{y_1, y_2\right\}$.
Since $\partial_D\left(V_3^{+-+}, y_i\right) \leq 2$ where $i=1,2,3$, then we have $V_3^{+-+} \rightarrow\left\{y_1, y_3\right\}$. \ Since $\partial_D\left(V_3^{-++}, y_i\right) \leq 2$ where $i=1,2,3$, then we have $V_3^{-++} \rightarrow\left\{y_2, y_3\right\}$. \ Since $\partial_D\left(x_3, V_3^{++-}\right) \leq 2, \partial_D\left(x_2, V_3^{+-+}\right) \leq 2$ and $\partial_D\left(x_1, V_3^{-++}\right) \leq 2$, then we have $y_3 \rightarrow V_3^{++-}, y_2 \rightarrow V_3^{+++}$and $y_1 \rightarrow V_3^{-++}$. \ Let $\left|V_3^{++-}\right|=q_1,\left|V_3^{+++}\right|=q_2,\left|V_3^{-++}\right|=q_3$. \ If $q_1 \geq 2$, then there exists $z_i, z_j \in V_3^{++-}$such that $\partial_D\left(z_i, z_j\right) \geq 3$, a contradiction. \ Hence $q_1 \leq 1$. \ If $q_2 \geq 2$, then there exists $z_i, z_j \in V_3^{+-+}$such that $\partial_D\left(z_i, z_j\right) \geq 3$, a contradiction. \ Hence $q_2 \leq 1$. \ If $q_3 \geq 2$, then there exists $z_i, z_j \in V_3^{-++}$such that $\partial_D\left(z_i, z_j\right) \geq 3$, a contradiction. \ Hence $q_3 \leq 1$. \ By Lemma 2.1, we also have $\left|V_3^{+++}\right| \leq 1$ and $\left|V_3^{---}\right| \leq 1$,\ thus $\left|V_3\right|=q \leq 1+1+1+1+1=5 \leq 6$.\\
(10) Case $(1,1,2)$.\\
Subcase1: $\left\{y_2, y_3\right\} \rightarrow\left\{x_1, x_2\right\} \rightarrow y_1, y_3 \rightarrow x_3 \rightarrow\left\{y_1, y_2\right\}$.

We know $y_1 \in V_2^{+++}, y_2 \in V_2^{--+}$and $y_3 \in V_2^{---}$. \ By Lemma 2.1, we can get $V_3^{+++}=V_3^{--+}=V_3^{---}=\varnothing$, thus $V_3=V_3^{++-} \cup V_3^{+-+} \cup V_3^{-++} \cup V_3^{+--} \cup V_3^{-+-}$. \ Take any $z \in V_3$, since $\partial_D\left(y_1, z\right) \leq 2$ and $\partial_D\left(z, y_3\right) \leq 2$, we have $y_1 \rightarrow z$ and $z \rightarrow y_3$, this means $y_1 \rightarrow V_3$ and $V_3 \rightarrow y_3 $. \ Since $\partial_D\left(V_3^{+-+}, y_i\right) \leq 2$ where $i=1,2,3$, then we have $V_3^{+-+} \rightarrow y_2$. \ Since $\partial_D\left(V_3^{-++}, y_i\right) \leq 2$ where $i=1,2,3$, then we have $V_3^{-++} \rightarrow y_2$. \ Let $\left|V_3^{++-}\right|=q_1,\left|V_3^{+-+}\right|=q_2,\left|V_3^{-++}\right|=q_3,\left|V_3^{+--}\right|=q_4$ and $\left|V_3^{-+-}\right|=q_5$. \ If $q_1 \geq 2$, then there exists $z_i, z_j \in V_3^{++-}$such that $\partial_D\left(z_i, z_j\right) \geq 3$, a contradiction. \ Hence $q_1 \leq 1$. \ If $q_2 \geq 2$, then there exists $z_i, z_j \in V_3^{+-+}$such that $\partial_D\left(z_i, z_j\right) \geq 3$, a contradiction. \ Hence $q_2 \leq 1$. If $q_3 \geq 2$, then there exists $z_i, z_j \in V_3^{-++}$such that $\partial_D\left(z_i, z_j\right) \geq 3$, a contradiction. \ Hence $q_3 \leq 1$. \ If $q_4 \geq 2$, then there exists $z_i, z_j \in V_3^{+--}$such that $\partial_D\left(z_i, z_j\right) \geq 3$, a contradiction. \ Hence $q_4 \leq 1$. \ If $q_5 \geq 2$, then there exists $z_i, z_j \in V_3^{-+-}$such that $\partial_D\left(z_i, z_j\right) \geq 3$, a contradiction. \ Hence $q_5 \leq 1$. \ Thus $\left|V_3\right|=q \leq 1+1+1+1+1=5 \leq 6$.\\
Subcase2: $\left\{y_2, y_3\right\} \rightarrow\left\{x_1, x_2\right\} \rightarrow y_1 \rightarrow x_3 \rightarrow\left\{y_2, y_3\right\}$.

We know $y_1 \in V_2^{++-}$and $y_2, y_3 \in V_2^{--+}$. \ By Lemma 2.1, we can get $V_3^{++-}=V_3^{--+}=\varnothing$, \ thus $V_3=V_3^{+++} \cup V_3^{-++} \cup V_3^{+-+} \cup V_3^{+--} \cup V_3^{-+-} \cup V_3^{---}$. \ Since $\partial_D\left(V_3^{+-+}, y_i\right) \leq 2$ where $i=1,2,3$, then we have $V_3^{+-+} \rightarrow\left\{y_2, y_3\right\}$. \ Since $\partial_D\left(V_3^{-++}, y_i\right) \leq 2$ where $i=1,2,3$, then we have $V_3^{-++} \rightarrow\left\{y_2, y_3\right\}$. \ Since $\partial_D\left(y_i, V_3^{+--}\right) \leq 2$ where $i=1,2,3$, then we have $y_1 \rightarrow V_3^{+--}$. \ Since $\partial_D\left(y_i, V_3^{-+-}\right) \leq 2$ where $i=1,2,3$, then we have $y_1 \rightarrow V_3^{-+-}$. \ Since $\partial_D\left(x_2, V_3^{+-+}\right) \leq 2$ and $\partial_D\left(x_1, V_3^{-++}\right) \leq 2$, we have $y_1 \rightarrow V_3^{+-+}$and $y_1 \rightarrow V_3^{-++}$. \ If $V_3^{+-+} \neq \varnothing$, then $\partial_D\left(V_3^{+-+}, x_3\right) \geq 3$, a contradiction. \ Hence $V_3^{+-+}=\varnothing$. \ Similarly, we have $V_3^{-++}=\varnothing$. \ Let $\left|V_3^{+--}\right|=q_1,\left|V_3^{-+-}\right|=q_2$ and $F_1=D\left[V_2^{--+} \cup V_3^{+--}\right], F_1=D\left[V_2^{--+} \cup V_3^{-+-}\right]$. \ Then $F_1$ and $F_2$ are respectively an orientation of $K\left(2, q_1\right)$ and $K\left(2, q_2\right)$ where $q_1, q_2 \leq q$. \ If $q_1 \geq 3$, then there exists $z_i, z_j \in V_3^{+--}$such that $\partial_{F_1}\left(z_i, z_j\right)=4$, so $\partial_D\left(z_i, z_j\right) \geq 3$, a contradiction. \ Hence $q_1 \leq 2$. \ If $q_2 \geq 3$, then there exists $z_i, z_j \in V_3^{-+-}$such that $\partial_{F_2}\left(z_i, z_j\right)=4$, so $\partial_D\left(z_i, z_j\right) \geq 3$, a contradiction. \ Hence $q_2 \leq 2$. \ By Lemma 2.1, we also have $\left|V_3^{+++}\right| \leq 1$ and $\left|V_3^{---}\right| \leq 1$, \ thus $\left|V_3\right|=q \leq 1+0+0+2+2+1=6 \leq 6$.\\
Subcase3: $\left\{y_2, y_3\right\} \rightarrow x_1 \rightarrow y_1,\left\{y_1, y_3\right\} \rightarrow x_2 \rightarrow y_2, y_3 \rightarrow x_3 \rightarrow\left\{y_1, y_2\right\}$.

We know $y_1 \in V_2^{+-+}, y_2 \in V_2^{-++}$and $y_3 \in V_2^{---}$. \ By Lemma 2.1, we can get $V_3^{+-+}=V_3^{-++}=V_3^{---}=\varnothing$, thus $V_3=V_3^{+++} \cup V_3^{++-} \cup V_3^{+--} \cup V_3^{-+-} \cup V_3^{--+}$. \ Take any $z \in V_3$, since $\partial_D\left(z, y_3\right) \leq 2$, then we have $z \rightarrow y_3$, this means $V_3 \rightarrow y_3$. \ Since $\partial_D\left(y_i, V_3^{+--}\right) \leq 2$, where $i=1,2,3$, then we have $y_1 \rightarrow V_3^{+--}$. \ Since $\partial_D\left(y_i, V_3^{-+-}\right) \leq 2$, where $i=1,2,3$, then we have $y_2 \rightarrow V_3^{-+-}$. \ Since $\partial_D\left(y_i, V_3^{--+}\right) \leq 2$, where $i=1,2,3$, then we have $\left\{y_1, y_2\right\} \rightarrow V_3^{--+}$. \ Let $\left|V_3^{++-}\right|=q_1,\left|V_3^{+--}\right|=q_2,\left|V_3^{-+-}\right|=q_3,\left|V_3^{--+}\right|=q_4$ and $F=D\left[V_2 \backslash\left\{y_3\right\} \cup V_3^{++-}\right]$. Then $F$ is an orientation of $K\left(2, q_1\right)$ where $q_1 \leq q$. \ If $q_1 \geq 3$, then there exists $z_i, z_j \in V_3^{++-}$such that $\partial_F\left(z_i, z_j\right)=4$, so $\partial_D\left(z_i, z_j\right) \geq 3$, a contradiction. \ Hence $q_1 \leq 2$. \ If $q_2 \geq 2$, then there exists $z_i, z_j \in V_3^{+--}$such that $\partial_D\left(z_i, z_j\right) \geq 3$, a contradiction. \ Hence $q_2 \leq 1$. \ If $q_3 \geq 2$, then there exists $z_i, z_j \in V_3^{-+-}$such that $\partial_D\left(z_i, z_j\right) \geq 3$, a contradiction. \ Hence $q_3 \leq 1$. If $q_4 \geq 2$, then there exists $z_i, z_j \in V_3^{--+}$such that $\partial_D\left(z_i, z_j\right) \geq 3$, a contradiction. \ Hence $q_4 \leq 1$. \ By Lemma 2.1, we also have $\left|V_3^{+++}\right| \leq 1$, \  thus $\left|V_3\right|=q \leq 1+2+1+1+1=6 \leq 6$.\\
Subcase 4: $\left\{y_2, y_3\right\} \rightarrow x_1 \rightarrow y_1,\left\{y_1, y_3\right\} \rightarrow x_2 \rightarrow y_2, y_1 \rightarrow x_3 \rightarrow\left\{y_2, y_3\right\}$.

We know $y_1 \in V_2^{+--}, y_2 \in V_2^{-++}$and $y_3 \in V_2^{--+}$. \ By Lemma 2.1, we can get
$V_3^{+--}=V_3^{-++}=V_3^{--+}=\varnothing$, thus $V_3=V_3^{+++} \cup V_3^{++-} \cup V_3^{+-+} \cup V_3^{-+-} \cup V_3^{---}$. \ Since $\partial_D\left(V_3^{++-}, y_i\right) \leq 2$ where $i=1,2,3$, we have $V_3^{++-} \rightarrow y_1$. \ Since $\partial_D\left(V_3^{+-+}, y_i\right) \leq 2$ where $i=1,2,3$, we have $V_3^{+-+} \rightarrow\left\{y_1, y_3\right\}$. \ Since $\partial_D\left(y_i, V_3^{-+-}\right) \leq 2$ where $i=1,2,3$, we have $y_2 \rightarrow V_3^{-+-}$. \ Since $\partial_D\left(x_2, V_3^{+-+}\right) \leq 2, \partial_D\left(x_1, V_3^{-+-}\right) \leq 2$ and $\partial_D\left(V_3^{-+-}, x_2\right) \leq 2$, then we have $y_2 \rightarrow V_3^{+-+}$and $y_1 \rightarrow V_3^{-+-} \rightarrow y_3$. \ Let $\left|V_3^{++-}\right|=q_1,\left|V_3^{+-+}\right|=q_2,\left|V_3^{-+-}\right|=q_3$ and $F=D\left[V_2 \backslash\left\{y_1\right\} \cup V_3^{++-}\right]$. \ Then $F$ is an orientation of $K\left(2, q_1\right)$ where $q_1 \leq q$. \ If $q_1 \geq 3$, then there exists $z_i, z_j \in V_3^{++-}$such that $\partial_F\left(z_i, z_j\right)=4$, so $\partial_D\left(z_i, z_j\right) \geq 3$, a contradiction. \ Hence $q_1 \leq 2$. \ If $q_2 \geq 2$, then there exists $z_i, z_j \in V_3^{+-+}$such that $\partial_D\left(z_i, z_j\right) \geq 3$, a contradiction. \ Hence $q_2 \leq 1$. If $q_3 \geq 2$, then there exists $z_i, z_j \in V_3^{-+-}$such that $\partial_D\left(z_i, z_j\right) \geq 3$, a contradiction. \ Hence $q_3 \leq 1$. \ By Lemma 2.1,\ we also have $\left|V_3^{+++}\right| \leq 1$ and $\left|V_3^{---}\right| \leq 1$, \ thus $\left|V_3\right|=q \leq 1+2+1+1+1=6 \leq 6$.   \\

In summary,\ it can be concluded that if $f(K(3,3, q))=2$, then $q \leq 6$. \  Since when $q \leq 6$, we have found an orientation of diameter 2 of $K(3,3, q)$ in Lemma 3.1. \ Therefore, $f(K(3,3, q))=2$ if and only if $q \leq 6$.  $\hfill\blacksquare$ \\

\section{ The oriented diameter of $K(3,4,q)$}\label{sec4}

\begin{lemma}
	For $4 \leq q \leq 11,  \  f(K(3,4, q))=2$.
	
\end{lemma}

\textit{\textbf{Proof.}} 
When $q=11$, let $V_3=V_3^{+++} \cup V_3^{++-} \cup V_3^{+-+} \cup V_3^{+--}$where $\left|V_3^{+++}\right|=1$, $\left|V_3^{++-}\right|=\left|V_3^{+-+}\right|=2$ and $\left|V_3^{+--}\right|=6$. Orient $K(3,4,11)$ as follows: $
V_2 \rightarrow x_1, y_4 \rightarrow x_2 \rightarrow V_2 \backslash\left\{y_4\right\}, y_1 \rightarrow x_3 \rightarrow V_2 \backslash\left\{y_1\right\}, V_3^{+++} \rightarrow V_2, V_3^{++-} \rightarrow\left\{y_1, y_4\right\},
$
$V_3^{+-+} \rightarrow\left\{y_1, y_4\right\}$. For the set $\left\{y_2, y_3\right\}$, let it form a directed cycle of length four with the vertices in $V_3^{++-}$and $V_3^{+-+}$respectively. The orientation between $V_2$ and $V_3^{+--}$is the same as the orientation of $K(4,6)$. Let $D_{11}$ be the resulting orientation, it is easy to verify that $\operatorname{diam}(D_{11})=2$.

When $q=10$, let $V_3=V_3^{+++} \cup V_3^{+-+} \cup V_3^{-++} \cup V_3^{+--} \cup V_3^{-+-} \cup V_3^{---}$where
$ \left|V_3^{+++}\right|=\left|V_3^{---}\right|=1,\left|V_3^{+-+}\right|=\left|V_3^{-++}\right|=\left|V_3^{+--}\right|=\left|V_3^{-+-}\right|=2 \text { and } V_3^{+++}=\left\{z_{+}\right\}, V_3^{---}=\left\{z_{-}\right\} , V_3^{+-+}=\left\{z_1, z_2\right\}, V_3^{-++}=\left\{z_3, z_4\right\}, V_3^{+--}=\left\{z_5, z_6\right\}, V_3^{-+-}=\left\{z_7, z_8\right\} .
$
Orient $K(3,4,10)$ as follows : $\left\{y_3, y_4\right\} \rightarrow\left\{x_1, x_2\right\} \rightarrow\left\{y_1, y_2\right\} \rightarrow x_3 \rightarrow\left\{y_3, y_4\right\}$,
$
V_3^{+++} \rightarrow V_2 \rightarrow V_3^{---},\left\{y_1, y_2\right\} \rightarrow V_3^{+--} \cup V_3^{-+-}, V_3^{+-+} \cup V_3^{-++} \rightarrow\left\{y_3, y_4\right\},
y_1 \rightarrow z_1 \rightarrow y_2 \rightarrow z_2 \rightarrow y_1, y_1 \rightarrow z_3 \rightarrow y_2 \rightarrow z_4 \rightarrow y_1, y_3 \rightarrow z_5 \rightarrow y_4 \rightarrow z_6 \rightarrow y_3, 
y_3 \rightarrow z_7 \rightarrow y_4 \rightarrow z_8 \rightarrow y_3 .
$
Let $D_{10}$ be the resulting orientation, it is easy to verify that $\operatorname{diam}(D_{10})=2$.

Let $D_{9}$ be the orientation obtained by deleting vertex $z_{-}$from the above orientation $D_{10}$, it is easy to verify $\operatorname{diam}\left(D_{9}\right)=2$, so $f(K(3,4,9))=2$.

Let $D_{8}$ be the orientation obtained by deleting vertex $z_{-}$and $z_{+}$ from the above orientation $D_{10}$, it is easy to verify $\operatorname{diam}\left(D_{8}\right)=2$, so $f(K(3,4,8))=2$.

Let $D_{7}$ be the orientation obtained by deleting vertex set $\left\{z_7, z_8\right\}$ and vertex $z_{-}$ from the above orientation $D_{10}$, it is easy to verify $\operatorname{diam}\left(D_{7}\right)=2$, so $f(K(3,4,7))=2$.

Let $D_{6}$ be the orientation obtained by deleting vertex set $\left\{z_7, z_8\right\}$ and vertex $z_{+}$, $z_{-}$ from the above orientation $D_{10}$, it is easy to verify $\operatorname{diam}\left(D_{6}\right)=2$, so
$ f(K(3,4,6))=2$.

Let $D_{5}$ be the orientation obtained by deleting vertex set $\left\{z_1, z_2\right\},\left\{z_7, z_8\right\}$ and vertex $z_{-}$ from the above orientation $D_{10}$, it is easy to verify $\operatorname{diam}\left(D_{5}\right)=2$, so
$f(K(3,4,5))=2$.

Let $D_{4}$ be the orientation obtained by deleting  vertex set $\left\{z_1, z_2\right\},\left\{z_7, z_8\right\}$ and vertex $z_{+}, z_{-}$ from the above orientation $D_{10}$, it is easy to verify $\operatorname{diam}\left(D_{4}\right)=2$, so $f(K(3,4,4))=2$.   $\hfill\blacksquare$ \\

\begin{theorem}
	Suppose $q \geq 4$, then
	$$
	f(K(3,4, q))= \begin{cases}2, & \text { if } q \leq 11 ; \\ 3, & \text { if } q>11. \end{cases}
	$$
	%\begin{eqnarray}\label{10}
	%....
	%\end{eqnarray}
\end{theorem}

\textit{\textbf{Proof.}}  When $4 \leq q \leq 11$,  we have shown in Lemma 4.1 that there exists an orientation of diameter 2 for $K(3,4, q)$. When $q>11$, we prove it by contradiction. \ Assuming $f(K(3,4, q))=2$ when $q>11$, then $K(3,4, q)$ has a strong orientation $D$ with diameter $\operatorname{diam}(D)=2$. \ Let $V_1=\left\{x_1, x_2, x_3\right\}, V_2=\left\{y_1, y_2, y_3, y_4\right\}$ and $V_3=\left\{z_1, z_2, \cdots z_q\right\}$ be the three parts of the vertex set of $K(3,4, q)$, and $i=\left|N_D^{+}\left(x_1\right) \cap V_2\right|, j=\left|N_D^{+}\left(x_2\right) \cap V_2\right|, k=\left|N_D^{+}\left(x_3\right) \cap V_2\right|$. \ For any $i, j, k \in\{0,1,2,3,4\}$, we have the following nineteen exhaustive cases Case(i,j,k)s to obtain the required contradiction.\\
(1) Case $(0,0,0). \ V_2 \rightarrow V_1$.

Take any $x \in V_1, z \in V_3$, since $\partial_D(x, z) \leq 2$, then we have $x \rightarrow z$, this means $V_1 \rightarrow V_3$. So $\quad \partial_D\left(x_1, x_2\right) \geq 3$, a contradiction.\\
(2) Case $(0,0,1). \  V_2 \rightarrow\left\{x_1, x_2\right\}, V_2 \backslash\left\{y_1\right\} \rightarrow x_3 \rightarrow y_1$.

Take any $x \in V_1 \backslash\left\{x_3\right\}$ and $z \in V_3$, since $\partial_D(x, z) \leq 2$, then we have $x \rightarrow z$, this means $V_1 \backslash\left\{x_3\right\} \rightarrow V_3$. So $\partial_D\left(x_1, x_2\right) \geq 3$, a contradiction.\\
(3) Case( $0,0,2). \  V_2 \rightarrow\left\{x_1, x_2\right\}, V_2 \backslash\left\{y_1, y_2\right\} \rightarrow x_3 \rightarrow\left\{y_1, y_2\right\}$.

Take any $x \in V_1 \backslash\left\{x_3\right\}$ and $z \in V_3$, since $\partial_D(x, z) \leq 2$, then we have $x \rightarrow z$, this means $V_1 \backslash\left\{x_3\right\} \rightarrow V_3$. So $\partial_D\left(x_1, x_2\right) \geq 3$, a contradiction.\\
(4) Case( $0,0,3). \  V_2 \rightarrow\left\{x_1, x_2\right\}, y_4 \rightarrow x_3 \rightarrow V_2 \backslash\left\{y_4\right\}$.

Take any $x \in V_1 \backslash\left\{x_3\right\}$ and $z \in V_3$, since $\partial_D(x, z) \leq 2$, then we have $x \rightarrow z$, this means $V_1 \backslash\left\{x_3\right\} \rightarrow V_3$. So $\partial_D\left(x_1, x_2\right) \geq 3$, a contradiction.\\
(5) Case( $0,0,4). \ V_2 \rightarrow\left\{x_1, x_2\right\}, x_3 \rightarrow V_2$.

Take any $x \in V_1 \backslash\left\{x_3\right\}$ and $z \in V_3$, since $\partial_D(x, z) \leq 2$, then we have $x \rightarrow z$, this means $V_1 \backslash\left\{x_3\right\} \rightarrow V_3$. So $\partial_D\left(x_1, x_2\right) \geq 3$, a contradiction.\\
(6) Case $(0,1,1)$.\\
Subcase 1: $V_2 \rightarrow x_1, V_2 \backslash\left\{y_1\right\} \rightarrow\left\{x_2, x_3\right\} \rightarrow y_1$.

Take any $y \in V_2 \backslash\left\{y_1\right\}$ and $z \in V_3$, since $\partial_D(z, y) \leq 2$, then we have $z \rightarrow y$, this means $V_3 \rightarrow V_2 \backslash\left\{y_1\right\}$. So $\partial_D\left(y_2, y_3\right) \geq 3$, a contradiction.\\
Subcase 2: $V_2 \rightarrow x_1, V_2 \backslash\left\{y_1\right\} \rightarrow x_2 \rightarrow y_1, V_2 \backslash\left\{y_2\right\} \rightarrow x_3 \rightarrow y_2$.

Take any $y \in V_2 \backslash\left\{y_1, y_2\right\}$ and $z \in V_3$, since $\partial_D(z, y) \leq 2$, then we have $z \rightarrow y$, this means $V_3 \rightarrow V_2 \backslash\left\{y_1, y_2\right\}$. So $\partial_D\left(y_3, y_4\right) \geq 3$, a contradiction.\\
(7) Case $(0,1,2)$.\\
Subcasse 1: $V_2 \rightarrow x_1, V_2 \backslash\left\{y_1\right\} \rightarrow x_2 \rightarrow y_1, V_2 \backslash\left\{y_1, y_2\right\} \rightarrow x_3 \rightarrow\left\{y_1, y_2\right\}$.

Take any $y \in V_2 \backslash\left\{y_1, y_2\right\}$ and $z \in V_3$, since $\partial_D(z, y) \leq 2$, then we have $z \rightarrow y$, this means $V_3 \rightarrow V_2 \backslash\left\{y_1, y_2\right\}$. So $\partial_D\left(y_3, y_4\right) \geq 3$, a contradiction.\\
Subcase 2: $V_2 \rightarrow x_1, V_2 \backslash\left\{y_1\right\} \rightarrow x_2 \rightarrow y_1, V_2 \backslash\left\{y_2, y_3\right\} \rightarrow x_3 \rightarrow\left\{y_2, y_3\right\}$.

We know $y_1 \in V_2^{-+-}, y_4 \in V_2^{---}$and $y_2, y_3 \in V_2^{--+}$. By Lemma 2.1, we can get $V_3^{-+-}=V_3^{--+}=V_3^{---}=\varnothing$. \ Take any $z \in V_3$, since $\partial_D\left(x_1, z\right) \leq 2$ and $\partial_D\left(z, y_4\right) \leq 2$, then we have $x_1 \rightarrow z$ and $z \rightarrow y_4$, this means $x_1 \rightarrow V_3$ and $V_3 \rightarrow y_4$. \ Thus $V_3=V_3^{+++} \cup V_3^{++-} \cup V_3^{+-+} \cup V_3^{+--}$. \ Since $\partial_D\left(V_3^{++-}, y_i\right) \leq 2$ where $i=1,2,3$, then we have $V_3^{++-} \rightarrow y_1$. \ Since $\partial_D\left(V_3^{+-+}, y_i\right) \leq 2$ where $i=1,2,3$, then we have $V_3^{+-+} \rightarrow\left\{y_2, y_3\right\}$. \ Since $\partial_D\left(x_2, V_3^{+-+}\right) \leq 2$, then we have $y_1 \rightarrow V_3^{+-+}$. \ Since $\partial_D\left(x_2, V_3^{+--}\right) \leq 2$, then we have $y_1 \rightarrow V_3^{+--}$. \ Let $\left|V_3^{++-}\right|=q_1,\left|V_3^{+-+}\right|=q_2,\left|V_3^{+--}\right|=q_3$ and $F_1=D\left[V_2^{--+} \cup V_3^{++-}\right], F_2=D\left[V_2^{--+} \cup V_3^{+--}\right]$. \ Then $F_1$ and $F_2$ are respextively an orientation of $K\left(2, q_1\right)$ and $K\left(2, q_3\right)$ where $q_1, q_3 \leq q$. \ If $q_1 \geq 3$, then there exists $z_i, z_j \in V_3^{++-}$such that $\partial_{F_1}\left(z_i, z_j\right)=4$, so $\partial_D\left(z_i, z_j\right) \geq 3$, a contradiction. \ Hence $q_1 \leq 2$. \ If $q_3 \geq 3$, then there exists $z_i, z_j \in V_3^{+--}$such that $\partial_{F_2}\left(z_i, z_j\right)=4$, so $\partial_D\left(z_i, z_j\right) \geq 3$, a contradiction. \ Hence $q_3 \leq 2$. \ If $q_2 \geq 2$, then there exists $z_i, z_j \in V_3^{+-+}$such that $\partial_D\left(z_i, z_j\right) \geq 3$, a contradiction. \ Hence $q_2 \leq 1$. \ By Lemma 2.1, we also have $\left|V_3^{+++}\right| \leq 1$, \ thus $\left|V_3\right|=q \leq 1+2+1+2=6 \leq 11$.\\
(8) Case $(0,1,3)$.\\
Subcase 1: $V_2 \rightarrow x_1, V_2 \backslash\left\{y_1\right\} \rightarrow x_2 \rightarrow y_1, y_4 \rightarrow x_3 \rightarrow V_2 \backslash\left\{y_4\right\}$.

We know $y_1 \in V_2^{-++}, y_4 \in V_2^{---}$and $y_2, y_3 \in V_2^{--+}$. By Lemma 2.1, we can get $V_3^{-++}=V_3^{--+}=V_3^{---}=\varnothing$. \ Take any $z \in V_3$, since $\partial_D\left(x_1, z\right) \leq 2$ and $\partial_D\left(z, y_4\right) \leq 2$, then we have $x_1 \rightarrow z$ and $z \rightarrow y_4$, this means $x_1 \rightarrow V_3$ and $V_3 \rightarrow y_4$. \ Thus $V_3=V_3^{+++} \cup V_3^{++-} \cup V_3^{+-+} \cup V_3^{+--}$. \ Since $\partial_D\left(V_3^{+-+}, y_i\right) \leq 2$ where $i=1,2,3$, then we have $V_3^{+-+} \rightarrow\left\{y_2, y_3\right\}$.  \ Since $\partial_D\left(x_2, V_3^{+-+}\right) \leq 2$, then we have $y_1 \rightarrow V_3^{+-+}$.  \ Since $\partial_D\left(x_2, V_3^{+--}\right) \leq 2$, then we have $y_1 \rightarrow V_3^{+--}$.  \ Let $\left|V_3^{++-}\right|=q_1,\left|V_3^{+-+}\right|=q_2,\left|V_3^{+--}\right|=q_3$ and $F_1=D\left[V_2 \backslash\left\{y_4\right\} \cup V_3^{++-}\right], F_2=D\left[V_2^{--+} \cup V_3^{+--}\right]$. \ Then $ F_1$ and $F_2$ are respextively an orientation of $K\left(3, q_1\right)$ and $K\left(2, q_3\right)$ where $q_1, q_3 \leq q$. \ If $q_1 \geq 4$, then there exists $z_i, z_j \in V_3^{++-}$such that $\partial_{F_1}\left(z_i, z_j\right)=4$, so $\partial_D\left(z_i, z_j\right) \geq 3$, a contradiction. \ Hence $q_1 \leq 3$. \ If $q_2 \geq 2$, then there exists $z_i, z_j \in V_3^{+-+}$such that $\partial_D\left(z_i, z_j\right) \geq 3$, a contradiction. \ Hence $q_2 \leq 1$. \ If $q_3 \geq 3$, then there exists $z_i, z_j \in V_3^{+--}$such that $\partial_{F_2}\left(z_i, z_j\right)=4$, so $\partial_D\left(z_i, z_j\right) \geq 3$, a contradiction. \ Hence $q_3 \leq 2$. \ By Lemma 2.1, we also have $\left|V_3^{+++}\right| \leq 1$, \ thus $\left|V_3\right|=q \leq 1+3+1+2=7 \leq 11$.\\
Subcase 2: $V_2 \rightarrow x_1, V_2 \backslash\left\{y_1\right\} \rightarrow x_2 \rightarrow y_1 \rightarrow x_3 \rightarrow V_2 \backslash\left\{y_1\right\}$.

We know $y_1 \in V_2^{-+-}$and $y_2, y_3, y_4 \in V_2^{--+}$. \ By Lemma 2.1, we can get $V_3^{-+-}=V_3^{--+}=\varnothing$. \ Take any $z \in V_3$, since $\partial_D\left(x_1, z\right) \leq 2$, then we have $x_1 \rightarrow z$, this means $x_1 \rightarrow V_3$. \ Thus $V_3=V_3^{+++} \cup V_3^{++-} \cup V_3^{+-+} \cup V_3^{+--}$. \ Since $\partial_D\left(V_3^{++-}, y_i\right) \leq 2$ where $i=1,2,3$, then we have $V_3^{++-} \rightarrow y_1$. \ Since $\partial_D\left(V_3^{+-+}, y_i\right) \leq 2$ where $i=1,2,3$, then we have $V_3^{+-+} \rightarrow\left\{y_2, y_3, y_4\right\}$. \ Since $\partial_D\left(x_2, V_3^{+-+}\right) \leq 2$, then we have $y_1 \rightarrow V_3^{+-+}$. \ Since $\partial_D\left(x_2, V_3^{+--}\right) \leq 2,  $then we have $y_1 \rightarrow V_3^{+--}$. \ If $V_3^{+-+} \neq \varnothing$, then $\partial_D\left(V_3^{+-+}, x_3\right) \geq 3$, a contradiction, so $V_3^{+-+}=\varnothing$. \ Let $\left|V_3^{++-}\right|=q_1,\left|V_3^{+--}\right|=q_2$ and $F_1=D\left[V_2^{--+} \cup V_3^{++-}\right], F_2=D\left[V_2^{--+} \cup V_3^{+--}\right]$. \ Then $F_1$ and $F_2$ are respextively an orientation of $K\left(3, q_1\right)$ and $K\left(3, q_2\right)$ where $q_1, q_2 \leq q$. \ If $q_1 \geq 4$, then there exists $z_i, z_j \in V_3^{++-}$such that $\partial_{F_1}\left(z_i, z_j\right)=4$, so $\partial_D\left(z_i, z_j\right) \geq 3$, a contradiction. \ Hence $q_1 \leq 3$. If $q_2 \geq 4$, then there exists $z_i, z_j \in V_3^{+--}$such that $\partial_{F_2}\left(z_i, z_j\right)=4$, so $\partial_D\left(z_i, z_j\right) \geq 3$, a contradiction. \ Hence $q_2 \leq 3$. \ By Lemma 2.1,  we also have $\left|V_3^{+++}\right| \leq 1$, thus $\left|V_3\right|=q \leq 1+3+0+3=7 \leq 11$.\\
(9) Case(0,1,4). $V_2 \rightarrow x_1, V_2 \backslash\left\{y_1\right\} \rightarrow x_2 \rightarrow y_1, x_3 \rightarrow V_2$.

We know $y_1 \in V_2^{-++}$and $y_2, y_3, y_4 \in V_2^{--+}$. \ By Lemma 2.1, we can get $V_3^{-++}=V_3^{--+}=\varnothing$. \ Take any $z \in V_3$, since $\partial_D\left(x_1, z\right) \leq 2$ and $\partial_D\left(z, x_3\right) \leq 2$, then we have $x_1 \rightarrow z$ and $z \rightarrow x_3$, this means $x_1 \rightarrow V_3$ and $V_3 \rightarrow x_3$. \ Thus $V_3=V_3^{++-} \cup V_3^{+--}$. \ Since $\partial_D\left(x_2, V_3^{+--}\right) \leq 2$, then we have $y_1 \rightarrow V_3^{+--}$. \ Let $\left|V_3^{++-}\right|=q_1,\left|V_3^{+--}\right|=q_2 $ and $F_1=D\left[V_2 \cup V_3^{++-}\right], F_2=D\left[V_2^{--+} \cup V_3^{+--}\right]$.  \ Then $ F_1$ and $F_2$ are respextively an orientation of $K\left(4, q_1\right)$ and $K\left(3, q_2\right)$ where $q_1, q_2 \leq q$.  \ If $q_1>\binom{4}{4 / 2}=6$, then there exists $z_i, z_j \in V_3^{++-}$such that $\partial_{F_1}\left(z_i, z_j\right)=4$, so $\partial_D\left(z_i, z_j\right) \geq 3$, a contradiction. \ Hence $q_1 \leq 6$. \ If $q_2 \geq 4$, then there exists $z_i, z_j \in V_3^{+--}$such that $\partial_{F_2}\left(z_i, z_j\right)=4$, so $\partial_D\left(z_i, z_j\right) \geq 3$, a contradiction. \ Hence $q_2 \leq 3$. \ Thus $\left|V_3\right|=q \leq 6+3=9 \leq 11$.\\
(10) Case $(0,2,2)$.\\
Subcase 1: $V_2 \rightarrow x_1, V_2 \backslash\left\{y_1, y_2\right\} \rightarrow\left\{x_2, x_3\right\} \rightarrow\left\{y_1, y_2\right\}$.

We know $y_1, y_2 \in V_2^{-++}$and $y_3, y_4 \in V_2^{---}$. \ By Lemma 2.1, we can get $V_3^{-++}=V_3^{---}=\varnothing$. \ Take any $z \in V_3$, since $\partial_D\left(x_1, z\right) \leq 2, \partial_D\left(z, y_3\right) \leq 2$ and $\partial_D\left(z, y_4\right) \leq 2$, then we have $x_1 \rightarrow z, z \rightarrow y_3$ and $z \rightarrow y_4$, this means $x_1 \rightarrow V_3$ and $V_3 \rightarrow\left\{y_3, y_4\right\}$. \ Thus $V_3=V_3^{+++} \cup V_3^{++-} \cup V_3^{+-+} \cup V_3^{+--}$. \ Let $\left|V_3^{++-}\right|=q_1,\left|V_3^{+-+}\right|=q_2,\left|V_3^{+--}\right|=q_3 $ and $F_1=D\left[V_2^{-++} \cup V_3^{++-}\right], F_2=D\left[V_2^{-++} \cup V_3^{+-+}\right], F_3=D\left[V_2^{-++} \cup V_3^{+--}\right]$. \ Then $F_1, F_2$ and are respextively an orientation of $K\left(2, q_1\right), K\left(2, q_2\right)$ and $K\left(2, q_3\right)$ where $q_1, q_2, q_3 \leq q$. \ If $q_1 \geq 3$, then there exists $z_i, z_j \in V_3^{++-}$such that $\partial_{F_1}\left(z_i, z_j\right)=4$, so $\partial_D\left(z_i, z_j\right) \geq 3$, a contradiction. \ Hence $q_1 \leq 2$. \ If $q_2 \geq 3$, then there exists $z_i, z_j \in V_3^{+-+}$such that $\partial_{F_2}\left(z_i, z_j\right)=4$, so $\partial_D\left(z_i, z_j\right) \geq 3$, a contradiction. \ Hence $q_2 \leq 2$. \ If $q_3 \geq 3$, then there exists $z_i, z_j \in V_3^{+--}$such that $\partial_{F_3}\left(z_i, z_j\right)=4$, so $\partial_D\left(z_i, z_j\right) \geq 3$, a contradiction. \ Hence $q_3 \leq 2$. \ By Lemma 2.1, we also have $\left|V_3^{+++}\right| \leq 1$, \ thus $\left|V_3\right|=q \leq 1+2+2+2=7 \leq 11$.\\
Subcase 2: $V_2 \rightarrow x_1, V_2 \backslash\left\{y_1, y_2\right\} \rightarrow x_2 \rightarrow\left\{y_1, y_2\right\}, V_2 \backslash\left\{y_2, y_3\right\} \rightarrow x_3 \rightarrow\left\{y_2, y_3\right\}$.

We know $y_1 \in V_2^{-+-}, y_2 \in V_2^{-++}, y_3 \in V_2^{--+}$and $y_4 \in V_2^{---}$. \ By Lemma 2.1, we can get $V_3^{-+-}=V_3^{-++}=V_3^{--+}=V_3^{---}=\varnothing$. \ Take any $z \in V_3$,  \ since $\partial_D\left(x_1, z\right) \leq 2$ and $\partial_D\left(z, y_4\right) \leq 2$, then we have $x_1 \rightarrow z$  and  $z \rightarrow y_4$, this means $x_1 \rightarrow V_3$  and  $V_3 \rightarrow y_4$. \ Thus $V_3=V_3^{+++} \cup V_3^{++-} \cup V_3^{+-+} \cup V_3^{+--}$. \ Since $\partial_D\left(V_3^{++-}, y_i\right) \leq 2$ where $i=1,2,3$, then we have $V_3^{++-} \rightarrow y_1$. \ Since $\partial_D\left(V_3^{+-+}, y_i\right) \leq 2$ where $i=1,2,3$, then we have $V_3^{+-+} \rightarrow y_3$. 
\ Let $\left|V_3^{++-}\right|=q_1,\left|V_3^{+-+}\right|=q_2,\left|V_3^{+--}\right|=q_3$ and $F_1=D\left[V_2 \backslash\left\{y_1, y_4\right\} \cup V_3^{++-}\right], F_2=D\left[V_2 \backslash\left\{y_3, y_4\right\} \cup V_3^{+-+}\right],  F_3=D\left[V_2 \backslash\left\{y_4\right\} \cup V_3^{+--}\right]$. \ Then $F_1, F_2$ and $F_3$ are respextively an orientation of $K\left(2, q_1\right), K\left(2, q_2\right)$ and $K\left(3, q_3\right)$  where   $q_1, q_2, q_3 \leq q$. \ If $q_1 \geq 3$, then there exists $z_i, z_j \in V_3^{++-}$such that $\partial_{F_1}\left(z_i, z_j\right)=4$, so $\partial_D\left(z_i, z_j\right) \geq 3$, a contradiction. \ Hence $q_1 \leq 2$. \ If $q_2 \geq 3$, then there exists $z_i, z_j \in V_3^{+-+}$such that $\partial_{F_2}\left(z_i, z_j\right)=4$, so $\partial_D\left(z_i, z_j\right) \geq 3$, a contradiction. \ Hence $q_2 \leq 2$. \ If $q_3 \geq 4$, then there exists $z_i, z_j \in V_3^{+--}$such that $\partial_{F_3}\left(z_i, z_j\right)=4$, so $\partial_D\left(z_i, z_j\right) \geq 3$, a contradiction. \ Hence $q_3 \leq 3$. \ By Lemma 2.1, we also have $\left|V_3^{+++}\right| \leq 1$, thus $\left|V_3\right|=q \leq 1+2+2+3=8 \leq 11$.\\
Subcase 3: $V_2 \rightarrow x_1, V_2 \backslash\left\{y_1, y_2\right\} \rightarrow x_2 \rightarrow\left\{y_1, y_2\right\}, V_2 \backslash\left\{y_3, y_4\right\} \rightarrow x_3 \rightarrow\left\{y_3, y_4\right\}$.

We know $y_1, y_2 \in V_2^{-+-}$and $y_3, y_4 \in V_2^{--+}$. \ By Lemma 2.1, we can get $V_3^{-+-}=V_3^{--+}=\varnothing$. \ Take any $z \in V_3$, since $\partial_D\left(x_1, z\right) \leq 2$, then we have $x_1 \rightarrow z$, this means $x_1 \rightarrow V_3$. Thus $V_3=V_3^{+++} \cup V_3^{++-} \cup V_3^{+-+} \cup V_3^{+--}$. \ Since $\partial_D\left(V_3^{++-}, y_i\right) \leq 2$ where $i=1,2,3$, then we have $V_3^{++-} \rightarrow\left\{y_1, y_2\right\}$. \ Since $\partial_D\left(V_3^{+-+}, y_i\right) \leq 2$ where $i=1,2,3$, then we have $V_3^{+-+} \rightarrow\left\{y_3, y_4\right\}$. \ Let $\left|V_3^{++-}\right|=q_1,\left|V_3^{+-+}\right|=q_2,\left|V_3^{+--}\right|=q_3$ and
$ F_1=D\left[V_2^{--+} \cup V_3^{++-}\right], F_2=D\left[V_2^{-+-} \cup V_3^{+-+}\right], F_3=D\left[V_2 \cup V_3^{+--}\right] .$  \ Then $F_1, F_2$ and $F_3$ are respextively an orientation of $K\left(2, q_1\right), K\left(2, q_2\right)$ and $K\left(4, q_3\right)$ where $q_1, q_2, q_3 \leq q$.  \ If $q_1 \geq 3$, then there exists $z_i, z_j \in V_3^{++-}$such that $\partial_{F_1}\left(z_i, z_j\right)=4$, so $\partial_D\left(z_i, z_j\right) \geq 3$, a contradiction. \ Hence $q_1 \leq 2$.  \ If $q_2 \geq 3$, then there exists $z_i, z_j \in V_3^{+-+}$such that $\partial_{F_2}\left(z_i, z_j\right)=4$, so $\partial_D\left(z_i, z_j\right) \geq 3$, a contradiction. \ Hence $q_2 \leq 2$.  \  If $q_3>\binom{4}{4 / 2}=6$, then there exists $z_i, z_j \in V_3^{+--}$such that $\partial_{F_3}\left(z_i, z_j\right)=4$, so $\partial_D\left(z_i, z_j\right) \geq 3$, a contradiction. \ Hence $q_3 \leq 6$. \ By Lemma 2.1, we also have $\left|V_3^{+++}\right| \leq 1$, \ thus $\left|V_3\right|=q \leq 1+2+2+6=11 \leq 11$.\\
(11) Case $(0,2,3)$.\\
Subcase 1: $V_2 \rightarrow x_1, V_2 \backslash\left\{y_1, y_2\right\} \rightarrow x_2 \rightarrow\left\{y_1, y_2\right\}, y_4 \rightarrow x_3 \rightarrow V_2 \backslash\left\{y_4\right\}$.

We know $y_1, y_2 \in V_2^{-++}, y_3 \in V_2^{--+}$and $y_4 \in V_2^{---}$. \ By Lemma 2.1, we can get $V_3^{-++}=V_3^{--+}=V_3^{---}=\varnothing$. \ Take any $z \in V_3$, since $\partial_D\left(x_1, z\right) \leq 2$ and $\partial_D\left(z, y_4\right) \leq 2$, then we have $x_1 \rightarrow z$ and $z \rightarrow y_4$, this means $x_1 \rightarrow V_3$ and $V_3 \rightarrow y_4$. \ Thus $V_3=V_3^{+++} \cup V_3^{++-} \cup V_3^{+-+} \cup V_3^{+--}$. \ Since $\partial_D\left(V_3^{+-+}, y_i\right) \leq 2$ where $i=1,2,3$, then we have $V_3^{+-+} \rightarrow y_3$. \ Let $\left|V_3^{++-}\right|=q_1,\left|V_3^{+-+}\right|=q_2,\left|V_3^{+--}\right|=q_3$ and $F_1=D\left[V_2 \backslash\left\{y_4\right\} \cup V_3^{++-}\right], F_2=D\left[V_2^{-++} \cup V_3^{+-+}\right], F_3=D\left[V_2 \backslash\left\{y_4\right\} \cup V_3^{+--}\right]$. \ Then $F_1, F_2$ and $F_3$ are respextively an orientation of $K\left(3, q_1\right), K\left(2, q_2\right)$ and $K\left(3, q_3\right)$ where $q_1, q_2, q_3\allowbreak  \leq q$.   
\ If $q_1 \geq 4$, then there exists $z_i, z_j \in V_3^{++-}$such that $\partial_{F_1}\left(z_i, z_j\right)=4$, so $\partial_D\left(z_i, z_j\right) \geq 3$, a contradiction. \ Hence $q_1 \leq 3$. \ If $q_2 \geq 3$, then there exists $z_i, z_j \in V_3^{+-+}$such that $\partial_{F_2}\left(z_i, z_j\right)=4$, so $\partial_D\left(z_i, z_j\right) \geq 3$, a contradiction. \ Hence $q_2 \leq 2$. \ If $q_3 \geq 4$, then there exists $z_i, z_j \in V_3^{+--}$such that $\partial_{F_3}\left(z_i, z_j\right)=4$, so $\partial_D\left(z_i, z_j\right) \geq 3$, a contradiction. \ Hence $q_3 \leq 3$. \ By Lemma 2.1, we also have $\left|V_3^{+++}\right| \leq 1$, \ thus $\left|V_3\right|=q \leq 1+3+2+3=9 \leq 11$.\\
Subcase 2: $V_2 \rightarrow x_1, V_2 \backslash\left\{y_1, y_2\right\} \rightarrow x_2 \rightarrow\left\{y_1, y_2\right\}, y_1 \rightarrow x_3 \rightarrow V_2 \backslash\left\{y_1\right\}$.

We know $y_1 \in V_2^{-+-}, y_2 \in V_2^{-++}$and $y_3, y_4 \in V_2^{--+}$. \ By Lemma 2.1, we can get $V_3^{-+-}=V_3^{-++}=V_3^{--+}=\varnothing$. \ Take any $z \in V_3$, since $\partial_D\left(x_1, z\right) \leq 2$, then we have $x_1 \rightarrow z$, this means $x_1 \rightarrow V_3$. \ Thus $V_3=V_3^{+++} \cup V_3^{++-} \cup V_3^{+-+} \cup V_3^{+--}$. \ Since $\partial_D\left(V_3^{++-}, y_i\right) \leq 2$ where $i=1,2,3$, then we have $V_3^{++-} \rightarrow y_1$. \ Since $\partial_D\left(V_3^{+-+}, y_i\right) \leq 2$ where $i=1,2,3$, then we have $V_3^{+-+} \rightarrow\left\{y_3, y_4\right\}$. \ Since $\partial_D\left(V_3^{+-+}, x_3\right) \leq 2$ and $\partial_D\left(x_2, V_3^{+-+}\right) \leq 2$, then we have $y_2 \rightarrow V_3^{+-+} \rightarrow y_1$. \ Let $\left|V_3^{++-}\right|=q_1,\left|V_3^{+-+}\right|=q_2,\left|V_3^{+--}\right|=q_3$ and $ F_1=D\left[V_2 \backslash\left\{y_1\right\} \cup V_3^{++-}\right], F_2=D\left[V_2 \cup V_3^{+--}\right]. \ $ Then $ F_1$ and $F_2$ are  respextively  an  orientation of $K\left(3, q_1\right)$ and $K\left(4, q_3\right)$ where $q_1, q_3 \leq q$.  
\ If $q_1 \geq 4$, then there exists $z_i, z_j \in V_3^{++-}$such that $\partial_{F_1}\left(z_i, z_j\right)=4$, so $\partial_D\left(z_i, z_j\right) \geq 3$, a contradiction. \ Hence $q_1 \leq 3$. If $q_2 \geq 2$, then there exists $z_i, z_j \in V_3^{+-+}$such that $\partial_D\left(z_i, z_j\right) \geq 3$, a contradiction. \ Hence $q_2 \leq 1$. \ If $q_3>\binom{4}{4 / 2}=6$, then there exists $z_i, z_j \in V_3^{+--}$such that $\partial_{F_2}\left(z_i, z_j\right)=4$, so $\partial_D\left(z_i, z_j\right) \geq 3$, a contradiction. \ Hence $q_3 \leq 6$. \ By Lemma 2.1, we also have $\left|V_3^{+++}\right| \leq 1$, \ thus $  \left|V_3\right|=q \leq 1+3+1+6=11 \leq 11 $.\\
(12) Case(0,2,4). $V_2 \rightarrow x_1, V_2 \backslash\left\{y_1, y_2\right\} \rightarrow x_2 \rightarrow\left\{y_1, y_2\right\}, x_3 \rightarrow V_2$.

We know $y_1, y_2 \in V_2^{-++}$and $y_3, y_4 \in V_2^{--+}$. \ By Lemma 2.1, we can get $V_3^{-++}=V_3^{--+}=\varnothing$. \ Take any $z \in V_3$, since $\partial_D\left(x_1, z\right) \leq 2$ and $\partial_D\left(z, x_3\right) \leq 2$, then we have $x_1 \rightarrow z$ and $z \rightarrow x_3$, this means $x_1 \rightarrow V_3$ and $V_3 \rightarrow x_3$. \ Thus $V_3=V_3^{++-} \cup V_3^{+--}$. \ Let $\left|V_3^{++-}\right|=q_1,\left|V_3^{+--}\right|=q_2$ and $F_1=D\left[V_2 \cup V_3^{++-}\right], F_2=D\left[V_2 \cup V_3^{+--}\right]$. \ Then $F_1$ and $F_2$ are respextively an orientation of $K\left(4, q_1\right)$ and $K\left(4, q_2\right)$ where $q_1, q_2 \leq q$. \ If $q_1>\binom{4}{4 / 2}=6$, then there exists $z_i, z_j \in V_3^{++-}$such that $\partial_{F_1}\left(z_i, z_j\right)=4$, so $\partial_D\left(z_i, z_j\right) \geq 3$, a contradiction. \ Hence $q_1 \leq 6$. \ If $q_2>\binom{4}{4 / 2}=6$, then there exists $z_i, z_j \in V_3^{+--}$such that $\partial_{F_2}\left(z_i, z_j\right)=4$, so $\partial_D\left(z_i, z_j\right) \geq 3$, a contradiction. \ Hence $q_2 \leq 6$. \ Since the orientation is unique when making $f(K(4,6))=2$, \ if $q_1=q_2=6$, then there exists $z_i \in V_3^{++-}, z_j \in V_3^{+--}$such that $\partial_D\left(z_i, z_j\right) \geq 3$, a contradiction. \ Hence we can get $q_1 \leq 5$ or $q_2 \leq 5$. \ The argument for these two cases are similar, so we may assume $q_2 \leq 5$. \ Thus $\left|V_3\right|=q \leq 6+5=11 \leq 11$.\\
(13) Case $(0,3,3)$.\\
Subcase 1: $V_2 \rightarrow x_1, y_4 \rightarrow\left\{x_2, x_3\right\} \rightarrow V_2 \backslash\left\{y_4\right\}$.

We know $y_1, y_2, y_3 \in V_2^{-++}$and $y_4 \in V_2^{---}$. \ By Lemma 2.1, we can get $V_3^{-++}=V_3^{---}=\varnothing$. \ Take any $z \in V_3$, since $\partial_D\left(x_1, z\right) \leq 2$ and $\partial_D\left(z, y_4\right) \leq 2$, then we have $x_1 \rightarrow z$ and $z \rightarrow y_4$, this means $x_1 \rightarrow V_3$ and $V_3 \rightarrow y_4$. \ Thus $V_3=V_3^{+++} \cup V_3^{++-} \cup V_3^{+-+} \cup V_3^{+--}$. \ Let $\left|V_3^{++-}\right|=q_1,\left|V_3^{+-+}\right|=q_2,\left|V_3^{+--}\right|=q_3$ and $F_1=D\left[V_2^{-++} \cup V_3^{++-}\right], F_2=D\left[V_2^{-++} \cup V_3^{+-+}\right], F_3=D\left[V_2^{-++} \cup V_3^{+--}\right]$. \ Then $F_1, F_2$ and $F_3$ are respextively an orientation of $K\left(3, q_1\right), K\left(3, q_2\right)$ and $K\left(3, q_3\right)$ where $q_1, q_2, q_3\allowbreak \leq q$.  \ If $q_1 \geq 4$, then there exists $z_i, z_j \in V_3^{++-}$such that $\partial_{F_1}\left(z_i, z_j\right)=4$, so $\partial_D\left(z_i, z_j\right) \geq 3$, a contradiction. \ Hence $q_1 \leq 3$. \ If $q_2 \geq 4$, then there exists $z_i, z_j \in V_3^{+-+}$such that $\partial_{F_2}\left(z_i, z_j\right)=4$, so $\partial_D\left(z_i, z_j\right) \geq 3$, a contradiction. \ Hence $q_2 \leq 3$. \ If $q_3 \geq 4$, then there exists $z_i, z_j \in V_3^{+--}$such that $\partial_{F_3}\left(z_i, z_j\right)=4$, so $\partial_D\left(z_i, z_j\right) \geq 3$, a contradiction. \ Hence $q_3 \leq 3$. \ By Lemma 2.1, we also have $\left|V_3^{+++}\right| \leq 1$, \ thus $\left|V_3\right|=q \leq 1+3+3+3=10 \leq 11$.\\
Subcase 2: $V_2 \rightarrow x_1, y_4 \rightarrow x_2 \rightarrow V_2 \backslash\left\{y_4\right\}, y_1 \rightarrow x_3 \rightarrow V_2 \backslash\left\{y_1\right\}$.

We know $y_1 \in V_2^{-+-}, y_2, y_3 \in V_2^{-++}$and $y_4 \in V_2^{--+}$. \ By Lemma 2.1, we can get $V_3^{-+-}=V_3^{-++}=V_3^{--+}=\varnothing$. \ Take any $z \in V_3$, since $\partial_D\left(x_1, z\right) \leq 2$, then we have $x_1 \rightarrow z$, this means $x_1 \rightarrow V_3$. \ Thus $V_3=V_3^{+++} \cup V_3^{++-} \cup V_3^{+-+} \cup V_3^{+--}$. \ Since $\partial_D\left(V_3^{++-}, y_i\right) \leq 2$ where $i=1,2,3$, then we have $V_3^{++-} \rightarrow y_1$. \ Since $\partial_D\left(V_3^{+-+}, y_i\right) \leq 2$ where $i=1,2,3$, then we have $V_3^{+-+} \rightarrow y_4$. \ Since $\partial_D\left(V_3^{++-}, x_2\right) \leq 2$, then we have $V_3^{++-} \rightarrow y_4$. \ Since $\partial_D\left(V_3^{+-+}, x_3\right) \leq 2$, then we have $V_3^{+-+} \rightarrow y_1$. \ Let $\left|V_3^{++-}\right|=q_1,\left|V_3^{+-+}\right|=q_2,\left|V_3^{+--}\right|=q_3$ and $F_1=D\left[V_2^{-++} \cup V_3^{++-}\right], F_2=D\left[V_2^{-++} \cup V_3^{+-+}\right], F_3=D\left[V_2 \cup V_3^{+--}\right]$. \ Then $F_1, F_2$ and $F_3$ are respextively an orientation of $K\left(2, q_1\right), K\left(2, q_2\right)$ and $K\left(4, q_3\right)$ where $q_1, q_2, q_3 \leq q$.  \ If $q_1 \geq 3$, then there exists $z_i, z_j \in V_3^{++-}$such that $\partial_{F_1}\left(z_i, z_j\right)=4$, so $\partial_D\left(z_i, z_j\right) \geq 3$, a contradiction. \ Hence $q_1 \leq 2$. \ If $q_2 \geq 3$, then there exists $z_i, z_j \in V_3^{+-+}$such that $\partial_{F_2}\left(z_i, z_j\right)=4$, so $\partial_D\left(z_i, z_j\right) \geq 3$, a contradiction. \ Hence $q_2 \leq 2$. \ If $q_3>\binom{4}{4 / 2}=6$, then there exists $z_i, z_j \in V_3^{+--}$such that $\partial_{F_3}\left(z_i, z_j\right)=4$, so $\partial_D\left(z_i, z_j\right) \geq 3$, a contradiction. \ Hence $q_3 \leq 6$. \ By Lemma 2.1, we also have $\left|V_3^{+++}\right| \leq 1$,  \ thus $\left|V_3\right|=q \leq 1+2+2+6=11 \leq 11$.\\
(14) Case $(1,1,1)$.\\
Subcase 1: $V_2 \backslash\left\{y_1\right\} \rightarrow V_1 \rightarrow y_1$.

Take any $y \in V_2 \backslash\left\{y_1\right\}$ and $z \in V_3$, since $\partial_D(z, y) \leq 2$, then we have $z \rightarrow y$, this means $V_3 \rightarrow V_2 \backslash\left\{y_1\right\}$. So $\partial_D\left(y_2, y_3\right) \geq 3$, a contradiction.\\
Subcase 2: $V_2 \backslash\left\{y_1\right\} \rightarrow\left\{x_1, x_2\right\} \rightarrow y_1, V_2 \backslash\left\{y_2\right\} \rightarrow x_3 \rightarrow y_2$.

Take any $y \in V_2 \backslash\left\{y_1, y_2\right\}$ and $z \in V_3$, since $\partial_D(z, y) \leq 2$, then we have $z \rightarrow y$, this means $V_3 \rightarrow V_2 \backslash\left\{y_1, y_2\right\}$. So $\partial_D\left(y_3, y_4\right) \geq 3$, a contradiction.\\
Subcase 3: $V_2 \backslash\left\{y_1\right\} \rightarrow x_1 \rightarrow y_1, V_2 \backslash\left\{y_2\right\} \rightarrow x_2 \rightarrow y_2, V_2 \backslash\left\{y_3\right\} \rightarrow x_3 \rightarrow y_3$.

We know $y_1 \in V_2^{+--}, y_2 \in V_2^{-+-}, y_3 \in V_2^{--+}$and $y_4 \in V_2^{---}$. \ By Lemma 2.1, we can get $V_3^{+--}=V_3^{-+-}=V_3^{--+}=V_3^{---}=\varnothing$. \ Thus $V_3=V_3^{+++} \cup V_3^{++-} \cup V_3^{+-+} \cup V_3^{-++}$. \ Take any $z \in V_3$, since $\partial_D\left(z, y_4\right) \leq 2$, then we have $z \rightarrow y_4$, this means $V_3 \rightarrow y_4$. \ Since $\partial_D\left(V_3^{++-}, y_i\right) \leq 2$ where $i=1,2,3$, then we have $V_3^{++-} \rightarrow\left\{y_1, y_2\right\}$. \ Since $\partial_D\left(V_3^{+-+}, y_i\right) \leq 2$ where $i=1,2,3$, then we have $V_3^{+-+} \rightarrow\left\{y_1, y_3\right\}$. \ Since $\partial_D\left(V_3^{-++}, y_i\right) \leq 2$ where $i=1,2,3$, then we have $V_3^{-++} \rightarrow\left\{y_2, y_3\right\}$. \ Since $\partial_D\left(x_3, V_3^{++-}\right) \leq 2$, then we have $y_3 \rightarrow V_3^{++-}$. \ Since $\partial_D\left(x_2, V_3^{+-+}\right) \leq 2$, then we have $y_2 \rightarrow V_3^{+-+}$. \ Since $\partial_D\left(x_1, V_3^{-++}\right) \leq 2$, then we have $y_1 \rightarrow V_3^{-++}$.  \ Let $\left|V_3^{++-}\right|=q_1,\left|V_3^{+-+}\right|=q_2$ and $\left|V_3^{-++}\right|=q_3$. \ If $q_1 \geq 2$, then there exists $z_i, z_j \in V_3^{++-}$such that $\partial_D\left(z_i, z_j\right) \geq 3$, a contradiction. \ Hence $q_1 \leq 1$. \ If $q_2 \geq 2$, then there exists $z_i, z_j \in V_3^{+-+}$such that $\partial_D\left(z_i, z_j\right) \geq 3$, a contradiction.  \ Hence $q_2 \leq 1$. \ If $q_3 \geq 2$, then there exists $z_i, z_j \in V_3^{-++}$such that $\partial_D\left(z_i, z_j\right) \geq 3$,  a contradiction. \ Hence $q_3 \leq 1$. \ By Lemma 2.1, we also have $\left|V_3^{+++}\right| \leq 1$, \ thus $\left|V_3\right|=q \leq 1+1+1+1=4 \leq 11$.\\
(15) Case $(1,1,2)$.\\
Subcase 1: $V_2 \backslash\left\{y_1\right\} \rightarrow\left\{x_1, x_2\right\} \rightarrow y_1, V_2 \backslash\left\{y_1, y_2\right\} \rightarrow x_3 \rightarrow\left\{y_1, y_2\right\}$.

Take any $y \in V_2 \backslash\left\{y_1, y_2\right\}$ and $z \in V_3$, since $\partial_D(z, y) \leq 2$, then we have $z \rightarrow y$, this means $V_3 \rightarrow V_2 \backslash\left\{y_1, y_2\right\}$. So $\partial_D\left(y_3, y_4\right) \geq 3$, a contradiction.\\
Subcase2: $V_2 \backslash\left\{y_1\right\} \rightarrow\left\{x_1, x_2\right\} \rightarrow y_1, V_2 \backslash\left\{y_2, y_3\right\} \rightarrow x_3 \rightarrow\left\{y_2, y_3\right\}$.

We know $y_1 \in V_2^{++-}, y_2, y_3 \in V_2^{--+}$and $y_4 \in V_2^{---}$. \ By Lemma 2.1, we can get $V_3^{++-}=V_3^{--+}=V_3^{---}=\varnothing$. \ Thus $V_3=V_3^{+++} \cup V_3^{+-+} \cup V_3^{-++} \cup V_3^{+--} \cup V_3^{-+-}$. \ Take any $z \in V_3$, since $\partial_D\left(z, y_4\right) \leq 2$, then we have $z \rightarrow y_4$, this means $V_3 \rightarrow y_4$. \ Since $\partial_D\left(V_3^{+-+}, y_i\right) \leq 2$ where $i=1,2,3$, then we have $V_3^{+-+} \rightarrow\left\{y_2, y_3\right\}$. \ Since $\partial_D\left(V_3^{-++}, y_i\right) \leq 2$ where $i=1,2,3$, then we have $V_3^{-++} \rightarrow\left\{y_2, y_3\right\}$. \ Since $\partial_D\left(x_2, V_3^{+--}\right) \leq 2$, then we have $y_1 \rightarrow V_3^{+--}$. \ Since $\partial_D\left(x_1, V_3^{-+-}\right) \leq 2$, then we have $y_1 \rightarrow V_3^{-+-}$.
\ Let $\left|V_3^{+-+}\right|=q_1,\left|V_3^{-++}\right|=q_2,\left|V_3^{+--}\right|=q_3,\left|V_3^{-+-}\right|=q_4$ and $F_1=D\left[V_2^{--+} \cup V_3^{+--}\right], F_2=D\left[V_2^{--+} \cup V_3^{-+-}\right]$. \ Then $ F_1 $ and $F_2$ are respextively an orientation of $K\left(2, q_3\right)$ and $K\left(2, q_4\right)$ where $q_3, q_4 \leq q$. \ If $q_1 \geq 2$, then there exists $z_i, z_j \in V_3^{+-+}$such that $\partial_D\left(z_i, z_j\right) \geq 3$, a contradiction. \ Hence $q_1 \leq 1$. \ If $q_2 \geq 2, \quad$ then there exists $z_i, z_j \in V_3^{-++}$such that $\partial_D\left(z_i, z_j\right) \geq 3,  \ $ a contradiction. \ Hence $q_2 \leq 1$. \ If $q_3 \geq 3$, then there exists $z_i, z_j \in V_3^{+--}$such that $\partial_{F_1}\left(z_i, z_j\right)=4$, so $\partial_D\left(z_i, z_j\right) \geq 3$, a contradiction. \ Hence $q_3 \leq 2$. If $q_4 \geq 3$, then there exists $z_i, z_j \in V_3^{-+-}$such that $\partial_{F_2}\left(z_i, z_j\right)=4$, so $\partial_D\left(z_i, z_j\right) \geq 3$, a contradiction. \ Hence $q_4 \leq 2.   \ $ By Lemma 2.1, we also have $\left|V_3^{+++}\right| \leq 1$, \ thus $ \left|V_3\right|=q \leq 1+1+1+2+2=7 \leq 11$.\\
Subcase 3: $V_2 \backslash\left\{y_1\right\} \rightarrow x_1 \rightarrow y_1, V_2 \backslash\left\{y_2\right\} \rightarrow x_2 \rightarrow y_2, V_2 \backslash\left\{y_1, y_2\right\} \rightarrow x_3 \rightarrow\left\{y_1, y_2\right\}$.

Take any $y \in V_2 \backslash\left\{y_1, y_2\right\}$ and $z \in V_3$, since $\partial_D(z, y) \leq 2$, then we have $z \rightarrow y$, this means $V_3 \rightarrow V_2 \backslash\left\{y_1, y_2\right\}$. So $\partial_D\left(y_3, y_4\right) \geq 3$, a contradiction.\\
Subcase 4: $V_2 \backslash\left\{y_1\right\} \rightarrow x_1 \rightarrow y_1, V_2 \backslash\left\{y_2\right\} \rightarrow x_2 \rightarrow y_2, V_2 \backslash\left\{y_2, y_3\right\} \rightarrow x_3 \rightarrow\left\{y_2, y_3\right\}$.

We know $y_1 \in V_2^{+--}, y_2 \in V_2^{-++}, y_3 \in V_2^{--+}$and $y_4 \in V_2^{---}$. \ By Lemma 2.1, we can get $V_3^{+--}=V_3^{-++}=V_3^{--+}=V_3^{---}=\varnothing$. Thus $V_3=V_3^{+++} \cup V_3^{++-} \cup V_3^{+-+} \cup V_3^{-+-}$. \ Take any $z \in V_3$, since $\partial_D\left(z, y_4\right) \leq 2$, then we have $z \rightarrow y_4$, this means $V_3 \rightarrow y_4$. \ Since $\partial_D\left(V_3^{++-}, y_i\right) \leq 2$ where $i=1,2,3$, then we have $V_3^{++-} \rightarrow y_1$. \ Since $\partial_D\left(V_3^{+-+}, y_i\right) \leq 2$ where $i=1,2,3$, then we have $V_3^{+-+} \rightarrow\left\{y_1, y_3\right\}$. \ Since $\partial_D\left(y_i, V_3^{-+-}\right) \leq 2$ where $i=1,2,3$, then we have $y_2 \rightarrow V_3^{-+-}$. \ Since $\partial_D\left(x_2, V_3^{+-+}\right) \leq 2$, then we have $y_2 \rightarrow V_3^{+-+}$. \ Since $\partial_D\left(x_1, V_3^{-+-}\right) \leq 2$, then we have $y_1 \rightarrow V_3^{-+-}$. \ Let $\left|V_3^{++-}\right|=q_1,\left|V_3^{+-+}\right|=q_2,\left|V_3^{-+-}\right|=q_3$ and $F=D\left[V_2 \backslash\left\{y_1, y_4\right\} \cup V_3^{++-}\right]$.
Then $F$ is an orientation of $K\left(2, q_1\right)$ where $q_1 \leq q$.
\ If $q_1 \geq 3$, then there exists $z_i, z_j \in V_3^{++-}$such that $\partial_F\left(z_i, z_j\right)=4$, so $\partial_D\left(z_i, z_j\right) \geq 3$, a contradiction. \ Hence $q_1 \leq 2$. \ If $q_2 \geq 2$, then there exists $z_i, z_j \in V_3^{+-+}$such that $\partial_D\left(z_i, z_j\right) \geq 3$, a contradiction. \ Hence $q_2 \leq 1$. \ If $q_3 \geq 2$, then there exists $z_i, z_j \in V_3^{-+-}$such that $\partial_D\left(z_i, z_j\right) \geq 3$, a contradiction. \ Hence $q_3 \leq 1$. \ By Lemma 2.1, we also have $\left|V_3^{+++}\right| \leq 1$, thus ||$V_3 \mid=q \leq 1+2+1+1=5 \leq 11$.\\
Subcase 5: $V_2 \backslash\left\{y_1\right\} \rightarrow x_1 \rightarrow y_1, V_2 \backslash\left\{y_2\right\} \rightarrow x_2 \rightarrow y_2, V_2 \backslash\left\{y_3, y_4\right\} \rightarrow x_3 \rightarrow\left\{y_3, y_4\right\}$.

We know $y_1 \in V_2^{+--}, y_2 \in V_2^{-+-}$and $y_3, y_4 \in V_2^{--+}$. \ By Lemma 2.1, we can get $V_3^{+--}=V_3^{-+-}=V_3^{--+}=\varnothing$. \ Thus $V_3=V_3^{+++} \cup V_3^{++-} \cup V_3^{+-+} \cup V_3^{-++} \cup V_3^{---}$. \ Since $\partial_D\left(V_3^{++-}, y_i\right) \leq 2$ where $i=1,2,3$, then we have $V_3^{++-} \rightarrow\left\{y_1, y_2\right\}$. \ Since $\partial_D\left(V_3^{+-+}, y_i\right) \leq 2$ where $i=1,2,3$, then we have $V_3^{+-+} \rightarrow\left\{y_1, y_3, y_4\right\}$. \ Since $\partial_D\left(V_3^{-++}, y_i\right) \leq 2$ where $i=1,2,3$, then we have $V_3^{-++} \rightarrow\left\{y_2, y_3, y_4\right\}$. \ Let $\left|V_3^{++-}\right|=q_1,\left|V_3^{+-+}\right|=q_2,\left|V_3^{-++}\right|=q_3 $  and $F=D\left[V_2^{--+} \cup V_3^{++-}\right]$. \ Then $F$ is an orientation of $K\left(2, q_1\right)$ where $q_1 \leq q$.  \ If $q_1 \geq 3$, then there exists $z_i, z_j \in V_3^{++-}$such that $\partial_F\left(z_i, z_j\right)=4$, so $\partial_D\left(z_i, z_j\right) \geq 3$, a contradiction. \ Hence $q_1 \leq 2$. \ If $q_2 \geq 2$, then there exists $z_i, z_j \in V_3^{+-+}$such that $\partial_D\left(z_i, z_j\right) \geq 3$, a contradiction. \ Hence $q_2 \leq 1$. If $q_3 \geq 2$, then there exists $z_i, z_j \in V_3^{-++}$such that $\partial_D\left(z_i, z_j\right) \geq 3$, a contradiction. \ Hence $q_3 \leq 1$. \ By Lemma 2.1, we also have $\left|V_3^{+++}\right| \leq 1$ and $\left|V_3^{---}\right| \leq 1$, \ thus $\left|V_3\right|=q \leq 1+2+1+1+1=6 \leq 11$.\\
(16) Case $(1,1,3)$.\\
Subcase 1: $V_2 \backslash\left\{y_1\right\} \rightarrow\left\{x_1, x_2\right\} \rightarrow y_1, y_4 \rightarrow x_3 \rightarrow V_2 \backslash\left\{y_4\right\}$.

We know $y_1 \in V_2^{+++}, y_2, y_3 \in V_2^{--+}$and $y_4 \in V_2^{---}$. \ By Lemma 2.1, we can get $V_3^{+++}=V_3^{--+}=V_3^{---}=\varnothing$. \ Thus $V_3=V_3^{++-} \cup V_3^{+-+} \cup V_3^{-++} \cup V_3^{+--} \cup V_3^{-+-}$. \ Take any $z \in V_3$, since $\partial_D\left(y_1, z\right) \leq 2$ and $\partial_D\left(z, y_4\right) \leq 2$, then we have $y_1 \rightarrow z$ and $z \rightarrow y_4$, this means $y_1 \rightarrow V_3$ and $V_3 \rightarrow y_4$. \ Since $\partial_D\left(V_3^{+-+}, y_i\right) \leq 2$ where $i=1,2,3$, then we have $V_3^{+-+} \rightarrow\left\{y_2, y_3\right\}$. \ Since $\partial_D\left(V_3^{-++}, y_i\right) \leq 2$ where $i=1,2,3$, then we have $V_3^{-++} \rightarrow\left\{y_2, y_3\right\}$. \ Let $\left|V_3^{++-}\right|=q_1,\left|V_3^{+-+}\right|=q_2,\left|V_3^{-++}\right|=q_3,\left|V_3^{+--}\right|=q_4,\left|V_3^{-+-}\right|=q_5   $ and $F_1=D\left[V_2^{--+} \cup V_3^{++-}\right], F_2=D\left[V_2^{--+} \cup V_3^{+--}\right], F_3=D\left[V_2^{--+} \cup V_3^{-+-}\right]$. \ Then $F_1, F_2$ and $F_3$ are respextively an orientation of $K\left(2, q_1\right), K\left(2, q_4\right)$ and $K\left(2, q_5\right)$ where $q_1, q_4, q_5 \leq q$. 
\ If $q_1 \geq 3$, then there exists $z_i, z_j \in V_3^{++-}$such that $\partial_{F_1}\left(z_i, z_j\right)=4$, so $\partial_D\left(z_i, z_j\right) \geq 3$, a contradiction. \ Hence $q_1 \leq 2$. \ If $q_2 \geq 2$, then there exists $z_i, z_j \in V_3^{+-+}$such that $\partial_D\left(z_i, z_j\right) \geq 3$, a contradiction. \ Hence $q_2 \leq 1$. \ If $q_3 \geq 2$, then there exists $z_i, z_j \in V_3^{-++}$such that $\partial_D\left(z_i, z_j\right) \geq 3$, a contradiction. \ Hence $q_3 \leq 1$. \ If $q_4 \geq 3$, then there exists $z_i, z_j \in V_3^{+--}$such that $\partial_{F_2}\left(z_i, z_j\right)=4$, so  $\partial_D\left(z_i, z_j\right) \geq 3$, a contradiction. \ Hence $q_4 \leq 2$. \ If $q_5 \geq 3$, then there exists $z_i, z_j \in V_3^{-+-}$such that $\partial_{F_3}\left(z_i, z_j\right)=4$, so $\partial_D\left(z_i, z_j\right) \geq 3$, a contradiction. \ Hence $q_5 \leq 2$. \ Thus $\left|V_3\right|=q \leq 2+1+1+2+2=8 \leq 11$.\\
Subcase 2: $V_2 \backslash\left\{y_1\right\} \rightarrow\left\{x_1, x_2\right\} \rightarrow y_1 \rightarrow x_3 \rightarrow V_2 \backslash\left\{y_1\right\}$.

We know $y_1 \in V_2^{++-}$and $y_2, y_3, y_4 \in V_2^{--+}$. \ By Lemma 2.1, we can get $V_3^{++-}=V_3^{--+}=\varnothing$. \ Thus $V_3=V_3^{+++} \cup V_3^{+-+} \cup V_3^{-++} \cup V_3^{+--} \cup V_3^{-+-} \cup V_3^{---}$. \ Since $\partial_D\left(V_3^{+-+}, y_i\right) \leq 2$ where $i=1,2,3$, then we have $V_3^{+-+} \rightarrow\left\{y_2, y_3, y_4\right\}$.
\ Since $\partial_D\left(V_3^{-++}, y_i\right) \leq 2$ where $i=1,2,3$, then we have $V_3^{-++} \rightarrow\left\{y_2, y_3, y_4\right\}$. \ Since $\partial_D\left(y_i, V_3^{+--}\right) \leq 2$ where $i=1,2,3$, then we have $y_1 \rightarrow V_3^{+--}$.
\ Since $\partial_D\left(y_i, V_3^{-+-}\right) \leq 2$ where $i=1,2,3$, then we have $y_1 \rightarrow V_3^{-+-}$. \ Let $\left|V_3^{+-+}\right|=q_1,\left|V_3^{-++}\right|=q_2,\left|V_3^{+--}\right|=q_3,\left|V_3^{-+-}\right|=q_4$ and $F_1=D\left[V_2^{--+} \cup V_3^{+--}\right], F_2=D\left[V_2^{--+} \cup V_3^{-+-}\right]. \ $ Then $ F_1$ and $F_2$ are respextively an orientation of $K\left(3, q_3\right)$ and $K\left(3, q_4\right)$ where $q_3, q_4 \leq q$.
\ If $q_1 \geq 2$, then there exists $z_i, z_j \in V_3^{+-+}$such that $\partial_D\left(z_i, z_j\right) \geq 3$, a contradiction. \ Hence $q_2 \leq 1$. \ If $q_2 \geq 2$, then there exists $z_i, z_j \in V_3^{-++}$such that $\partial_D\left(z_i, z_j\right) \geq 3$, a contradiction. \ Hence $q_2 \leq 1$. If $q_3 \geq 4$, then there exists $z_i, z_j \in V_3^{+--}$such that $\partial_{F_1}\left(z_i, z_j\right)=4$, so $\partial_D\left(z_i, z_j\right) \geq 3$, a contradiction. \ Hence $q_3 \leq 3$.\ If $q_4 \geq 4$, then there exists $z_i, z_j \in V_3^{-+-}$such that $\partial_{F_2}\left(z_i, z_j\right)=4$, so $\partial_D\left(z_i, z_j\right) \geq 3$, a contradiction. \ Hence $q_4 \leq 3$. \ By Lemma 2.1, we also have $\left|V_3^{+++}\right| \leq 1$ and $\left|V_3^{---}\right| \leq 1$, \ thus $\left|V_3\right|=q \leq 1+1+1+3+3+1=10 \leq 11$.\\
Subcase 3: $V_2 \backslash\left\{y_1\right\} \rightarrow x_1 \rightarrow y_1, V_2 \backslash\left\{y_2\right\} \rightarrow x_2 \rightarrow y_2, y_4 \rightarrow x_3 \rightarrow V_2 \backslash\left\{y_4\right\}$.

We know $y_1 \in V_2^{+-+}, y_2 \in V_2^{-++}, y_3 \in V_2^{--+}$and $y_4 \in V_2^{---}$.  \ By Lemma 2.1, we can get $V_3^{+-+}=V_3^{-++}=V_3^{--+}=V_3^{---}=\varnothing$. \ Thus $V_3=V_3^{+++} \cup V_3^{++-} \cup V_3^{+--} \cup V_3^{-+-}$. \ Take any $z \in V_3$, since $\partial_D\left(z, y_4\right) \leq 2$, then we have $z \rightarrow y_4$, this means $V_3 \rightarrow y_4$. \ Since $\partial_D\left(y_i, V_3^{+--}\right) \leq 2$ where $i=1,2,3$, then we have $y_1 \rightarrow V_3^{+--}$. \ Since $\partial_D\left(y_i, V_3^{-+-}\right) \leq 2$ where $i=1,2,3$, then we have $y_2 \rightarrow V_3^{-+-}$. \ Let $\left|V_3^{++-}\right|=q_1,\left|V_3^{+--}\right|=q_2,\left|V_3^{-+-}\right|=q_3$ and $
F_1=D\left[V_2 \backslash\left\{y_4\right\} \cup V_3^{++-}\right], F_2=D\left[V_2 \backslash\left\{y_1, y_4\right\} \cup V_3^{+--}\right], F_3=D\left[V_2 \backslash\left\{y_2, y_4\right\} \cup V_3^{-+-}\right].$ \ Then $F_1, F_2$ and $F_3$ are respextively an orientation of $K\left(3, q_1\right), K\left(2, q_2\right)$ and $K\left(2, q_3\right)$ where $q_1, q_2, q_3 \leq q$. \ If $q_1 \geq 4$, then there exists $z_i, z_j \in V_3^{++-}$such that $\partial_{F_1}\left(z_i, z_j\right)=4$, so $\partial_D\left(z_i, z_j\right) \geq 3$, a contradiction. \ Hence $q_1 \leq 3$. \ If $q_2 \geq 3$, then there exists $z_i, z_j \in V_3^{+--}$such that $\partial_{F_2}\left(z_i, z_j\right)=4$, so $\partial_D\left(z_i, z_j\right) \geq 3$, a contradiction. \ Hence $q_2 \leq 2$. \ If $q_3 \geq 3$, then there exists $z_i, z_j \in V_3^{-+-}$such that $\partial_{F_3}\left(z_i, z_j\right)=4$, so $\partial_D\left(z_i, z_j\right) \geq 3$, a contradiction. \ Hence $q_3 \leq 2$. \ By Lemma 2.1, we also have $\left|V_3^{+++}\right| \leq 1,  \ $ thus $\left|V_3\right|=q \leq 1+3+2+2=8 \leq 11$.\\
Subcase 4: $V_2 \backslash\left\{y_1\right\} \rightarrow x_1 \rightarrow y_1, V_2 \backslash\left\{y_2\right\} \rightarrow x_2 \rightarrow y_2, y_1 \rightarrow x_3 \rightarrow V_2 \backslash\left\{y_1\right\}$.

We know $y_1 \in V_2^{+--}, y_2 \in V_2^{-++}$and $y_3, y_4 \in V_2^{--+}$. \ By Lemma 2.1, we can get $V_3^{+--}=V_3^{-++}=V_3^{--+}=\varnothing$. \ Thus $V_3=V_3^{+++} \cup V_3^{++-} \cup V_3^{+-+} \cup V_3^{-+-} \cup V_3^{---}$. \ Since $\partial_D\left(V_3^{++-}, y_i\right) \leq 2$ where $i=1,2,3$, then we have $V_3^{++-} \rightarrow y_1$. \ Since $\partial_D\left(V_3^{+-+}, y_i\right) \leq 2$ where $i=1,2,3$, then we have $V_3^{+-+} \rightarrow\left\{y_1, y_3, y_4\right\}$. \ Since $\partial_D\left(y_i, V_3^{-+-}\right) \leq 2$ where $i=1,2,3$, then we have $y_2 \rightarrow V_3^{-+-}$.  \ Let $\left|V_3^{++-}\right|=q_1,\left|V_3^{+-+}\right|=q_2,\left|V_3^{-+-}\right|=q_3$ and $F_1=D\left[V_2 \backslash\left\{y_1\right\} \cup V_3^{++-}\right], F_2=D\left[V_2 \backslash\left\{y_2\right\} \cup V_3^{-+-}\right]$. \ Then $F_1$ and $F_2$ are respextively an orientation of $K\left(3, q_1\right)$ and $K\left(3, q_3\right)$ where $q_1, q_3 \leq q$. \ If $q_1 \geq 4$, then there exists $z_i, z_j \in V_3^{++-}$such that $\partial_{F_1}\left(z_i, z_j\right)=4$, so $\partial_D\left(z_i, z_j\right) \geq 3$, a contradiction. \ Hence $q_1 \leq 3$. \ If $q_2 \geq 2$, then there exists $z_i, z_j \in V_3^{+-+}$such that $\partial_D\left(z_i, z_j\right) \geq 3$, a contradiction. \ Hence $q_2 \leq 1$. \ If $q_3 \geq 4$, then there exists $z_i, z_j \in V_3^{-+-}$such that $\partial_{F_2}\left(z_i, z_j\right)=4$, so $\partial_D\left(z_i, z_j\right) \geq 3$, a contradiction. \ Hence $q_3 \leq 3$. \ By Lemma 2.1, we also have $\left|V_3^{+++}\right| \leq 1$ and $\left|V_3^{---}\right| \leq 1$, \ thus $\left|V_3\right|=q \leq 1+3+1+3+1=9 \leq 11$.\\
(17) Case( $1,2,2)$.\\
Subcase 1: $V_2 \backslash\left\{y_1\right\} \rightarrow x_1 \rightarrow y_1, V_2 \backslash\left\{y_1, y_2\right\} \rightarrow\left\{x_2, x_3\right\} \rightarrow\left\{y_1, y_2\right\}$.

Take any $z \in V_3$, since $\partial_D\left(z, y_3\right) \leq 2$ and $\partial_D\left(z, y_4\right) \leq 2$, then we have $z \rightarrow y_3$ and $z \rightarrow y_4$, this means $V_3 \rightarrow\left\{y_3, y_4\right\}$. So $\partial_D\left(y_3, y_4\right) \geq 3$, a contradiction.\\
Subcase 2: $V_2 \backslash\left\{y_1\right\} \rightarrow x_1 \rightarrow y_1,\left\{y_3, y_4\right\} \rightarrow x_2 \rightarrow\left\{y_1, y_2\right\},\left\{y_2, y_4\right\} \rightarrow x_3 \rightarrow\left\{y_1, y_3\right\}$.

We know $y_1 \in V_2^{+++}, y_2 \in V_2^{-+-}, y_3 \in V_2^{--+}$and $y_4 \in V_2^{---}$.  \ By Lemma 2.1, we can get $V_3^{+++}=V_3^{-+-}=V_3^{--+}=V_3^{---}=\varnothing$. \ Thus $V_3=V_3^{++-} \cup V_3^{+-+} \cup V_3^{-++} \cup V_3^{+--}$. \ Take any $z \in V_3$, since $\partial_D\left(y_1, z\right) \leq 2$ and $\partial_D\left(z, y_4\right) \leq 2$, then we have $y_1 \rightarrow z$ and $z \rightarrow y_4$, this means $y_1 \rightarrow V_3$ and $V_3 \rightarrow y_4$.  \ Since $\partial_D\left(V_3^{++-}, y_i\right) \leq 2$ where $i=1,2,3$, then we have $V_3^{++-} \rightarrow y_2$. \ Since $\partial_D\left(V_3^{+-+}, y_i\right) \leq 2$ where $i=1,2,3$, then we have $V_3^{+-+} \rightarrow y_3$. \ Since $\partial_D\left(V_3^{-++}, y_i\right) \leq 2$ where $i=1,2,3$, then we have $V_3^{-++} \rightarrow\left\{y_2, y_3\right\}$. \ Let $\left|V_3^{++-}\right|=q_1,\left|V_3^{+-+}\right|=q_2,\left|V_3^{-++}\right|=q_3,\left|V_3^{+--}\right|=q_4$ and $F=D\left[V_2 \backslash\left\{y_1, y_4\right\} \cup V_3^{+--}\right]$. \ Then $F$ is an orientation of $K\left(2, q_4\right)$ where $q_4 \leq q$. \ If $q_1 \geq 2$, then there exists $z_i, z_j \in V_3^{++-}$such that $\partial_D\left(z_i, z_j\right) \geq 3,\ $ a contradiction. \ Hence $q_1 \leq 1$. \ If $q_2 \geq 2$, then there exists $z_i, z_j \in V_3^{+-+}$such that $\partial_D\left(z_i, z_j\right) \geq 3$, a contradiction. \ Hence $q_2 \leq 1$. \ If $q_3 \geq 2$, then there exists $z_i, z_j \in V_3^{-++}$such that $\partial_D\left(z_i, z_j\right) \geq 3$, a contradiction. \ Hence $q_3 \leq 1$. \ If $q_4 \geq 3$, then there exists $z_i, z_j \in V_3^{+--}$such that $\partial_F\left(z_i, z_j\right)=4$, so $\partial_D\left(z_i, z_j\right) \geq 3$, a contradiction. \ Hence $q_4 \leq 2$. \ Thus $\left|V_3\right|=q \leq 1+1+1+2=5 \leq 11$.\\
Subcase 3: $V_2 \backslash\left\{y_1\right\} \rightarrow x_1 \rightarrow y_1,\left\{y_3, y_4\right\} \rightarrow x_2 \rightarrow\left\{y_1, y_2\right\},\left\{y_1, y_4\right\} \rightarrow x_3 \rightarrow\left\{y_2, y_3\right\}$.

We know $y_1 \in V_2^{++-}, y_2 \in V_2^{-++}, y_3 \in V_2^{--+}$and $y_4 \in V_2^{---}$.  \ By Lemma 2.1, we can get $V_3^{++-}=V_3^{-++}=V_3^{--+}=V_3^{---}=\varnothing$. \ Thus $V_3=V_3^{+++} \cup V_3^{+-+} \cup V_3^{+--} \cup V_3^{-+-}$. \ Take any $z \in V_3$, since $\partial_D\left(z, y_4\right) \leq 2$, then we have $z \rightarrow y_4$, this means $V_3 \rightarrow y_4$. \ Since $\partial_D\left(V_3^{+-+}, y_i\right) \leq 2$ where $i=1,2,3$, then we have $V_3^{+-+} \rightarrow y_3$. \ Since $\partial_D\left(y_i, V_3^{+--}\right) \leq 2$ where $i=1,2,3$, then we have $y_1 \rightarrow V_3^{+--}$. \ Since $\partial_D\left(y_i, V_3^{-+-}\right) \leq 2$ where $i=1,2,3$, then we have $\left\{y_1, y_2\right\} \rightarrow V_3^{-+-}$. \ Let $\left|V_3^{+-+}\right|=q_1,\left|V_3^{+--}\right|=q_2,\left|V_3^{-+-}\right|=q_3$ and $\mid F_1=D\left[V_2 \backslash\left\{y_3, y_4\right\} \cup V_3^{+-+}\right], F_2=D\left[V_2 \backslash\left\{y_1, y_4\right\} \cup V_3^{+--}\right]$. \ Then $F_1$ and $F_2$ are respextively an orientation of $K\left(2, q_1\right)$ and $K\left(2, q_2\right)$ where $q_1, q_2 \leq q$.  \ If $q_1 \geq 3$, then there exists $z_i, z_j \in V_3^{+-+}$such that $\partial_{F_1}\left(z_i, z_j\right)=4$, so $\partial_D\left(z_i, z_j\right) \geq 3$, a contradiction. \ Hence $q_1 \leq 2$. \ If $q_2 \geq 3$, then there exists $z_i, z_j \in V_3^{+--}$such that $\partial_{F_2}\left(z_i, z_j\right)=4$, so $\partial_D\left(z_i, z_j\right) \geq 3$, a contradiction. \ Hence $q_2 \leq 2$. \ If $q_3 \geq 2$, then there exists $z_i, z_j \in V_3^{-+-}$such that $\partial_D\left(z_i, z_j\right) \geq 3$, a contradiction. \ Hence $q_3 \leq 1$. \ By Lemma 2.1, we also have $\left|V_3^{+++}\right| \leq 1$, \ thus $\left|V_3\right|=q \leq 1+2+2+1=6 \leq 11$.\\
Subcase 4: $V_2 \backslash\left\{y_1\right\} \rightarrow x_1 \rightarrow y_1,\left\{y_3, y_4\right\} \rightarrow x_2 \rightarrow\left\{y_1, y_2\right\} \rightarrow x_3 \rightarrow\left\{y_3, y_4\right\}$.

We know $y_1 \in V_2^{++-}, y_2 \in V_2^{-+-}$and $y_3, y_4 \in V_2^{--+}$. \ By Lemma 2.1, we can get $V_3^{++-}=V_3^{-+-}=V_3^{--+}=\varnothing$. \ Thus $V_3=V_3^{+++} \cup V_3^{+-+} \cup V_3^{-++} \cup V_3^{+--} \cup V_3^{---}$. \ Since $\partial_D\left(V_3^{+-+}, y_i\right) \leq 2$ where $i=1,2,3$, then we have $V_3^{+-+} \rightarrow\left\{y_3, y_4\right\}$. \ Since $\partial_D\left(V_3^{-++}, y_i\right) \leq 2$ where $i=1,2,3$, then we have $V_3^{-++} \rightarrow\left\{y_2, y_3, y_4\right\}$. \ Since $\partial_D\left(y_i, V_3^{+--}\right) \leq 2$ where $i=1,2,3$, then we have $y_1 \rightarrow V_3^{+--}$. \ Let $\left|V_3^{+-+}\right|=q_1,\left|V_3^{-++}\right|=q_2,\left|V_3^{+--}\right|=q_3$ and $F_1=D\left[V_2 \backslash\left\{y_3, y_4\right\} \cup V_3^{+-+}\right], F_2=D\left[V_2 \backslash\left\{y_1\right\} \cup V_3^{+--}\right]$. \ Then $F_1$ and $F_2$ are respextively an orientation of $K\left(2, q_1\right)$ and $K\left(3, q_3\right)$ where $q_1, q_3 \leq q$.  \ If $q_1 \geq 3$, then there exists $z_i, z_j \in V_3^{+-+}$such that $\partial_{F_1}\left(z_i, z_j\right)=4$, so $\partial_D\left(z_i, z_j\right) \geq 3$, a contradiction. \ Hence $q_1 \leq 2$. \ If $q_2 \geq 2$, then there exists $z_i, z_j \in V_3^{-++}$such that $\partial_D\left(z_i, z_j\right) \geq 3$, a contradiction. \ Hence $q_2 \leq 1$. \ If $q_3 \geq 4$, then there exists $ z_i, z_j \in V_3^{+--}$such that $\partial_{F_2}\left(z_i, z_j\right)=4$, so $\partial_D\left(z_i, z_j\right) \geq 3$, a contradiction. \ Hence $q_3 \leq 3$. \ By Lemma 2.1, we also have $\left|V_3^{+++}\right| \leq 1$ and $\left|V_3^{---}\right| \leq 1$, thus $\left|V_3\right|=q \leq 1+2+1+3+1=8 \leq 11$.\\
Subcase 5: $V_2 \backslash\left\{y_1\right\} \rightarrow x_1 \rightarrow y_1, V_2 \backslash\left\{y_2, y_3\right\} \rightarrow\left\{x_2, x_3\right\} \rightarrow\left\{y_2, y_3\right\}$.

We know $y_1 \in V_2^{+--}, y_2, y_3 \in V_2^{-++}$and $y_4 \in V_2^{---}$. \ By Lemma 2.1, we can get  $V_3^{+--}=V_3^{-++}=V_3^{---}=\varnothing$. \ Thus $V_3=V_3^{+++} \cup V_3^{++-} \cup V_3^{+-+} \cup V_3^{-+-} \cup V_3^{--+}$.  \ Take any $z \in V_3$, since $\partial_D\left(z, y_4\right) \leq 2$, then we have $z \rightarrow y_4$, this means $V_3 \rightarrow y_4$. \ Since $\partial_D\left(V_3^{++-}, y_i\right) \leq 2$ where $i=1,2,3$, then we have $V_3^{++-} \rightarrow y_1$. \ Since $\partial_D\left(V_3^{+-+}, y_i\right) \leq 2$ where $i=1,2,3$, then we have $V_3^{+-+} \rightarrow y_1$. \ Since $\partial_D\left(y_i, V_3^{-+-}\right) \leq 2$ where $i=1,2,3$, then we have $\left\{y_2, y_3\right\} \rightarrow V_3^{-+-}$. \ Since $\partial_D\left(y_i, V_3^{--+}\right) \leq 2$ where $i=1,2,3$, then we have $\left\{y_2, y_3\right\} \rightarrow V_3^{--+}$. \ Let $\left|V_3^{++-}\right|=q_1,\left|V_3^{+-+}\right|=q_2,\left|V_3^{-+-}\right|=q_3,\left|V_3^{--+}\right|=q_4$ and $F_1=D\left[V_2^{-++} \cup V_3^{++-}\right], F_2=D\left[V_2^{-++} \cup V_3^{+-+}\right]$. \ Then $F_1$ and $F_2$ are respextively an orientation of $K\left(2, q_1\right)$ and $K\left(2, q_2\right)$ where $q_1, q_2 \leq q$. \ If $q_1 \geq 3$, then there exists $z_i, z_j \in V_3^{++-}$such that $\partial_{F_1}\left(z_i, z_j\right)=4$, so $\partial_D\left(z_i, z_j\right) \geq 3$, a contradiction. \ Hence $q_1 \leq 2$. \ If $q_2 \geq 3$, then there exists $z_i, z_j \in V_3^{+-+}$such that $\partial_{F_2}\left(z_i, z_j\right)=4$, so $\partial_D\left(z_i, z_j\right) \geq 3$, a contradiction. \ Hence $q_2 \leq 2$. \ If $q_3 \geq 2$, then there exists $z_i, z_j \in V_3^{-+-}$such that $\partial_D\left(z_i, z_j\right) \geq 3$, a contradiction. \ Hence $q_3 \leq 1$. If $q_4 \geq 2$, then there exists $z_i, z_j \in V_3^{--+}$such that $\partial_D\left(z_i, z_j\right) \geq 3$, a contradiction. \ Hence $q_4 \leq 1$. \ By Lemma 2.1, we also have $\left|V_3^{+++}\right| \leq 1$, thus $\left|V_3\right|=q \leq 1+2+2+1+1=7 \leq 11$.\\
Subcase 6: $V_2 \backslash\left\{y_1\right\} \rightarrow x_1 \rightarrow y_1,\left\{y_1, y_4\right\} \rightarrow x_2 \rightarrow\left\{y_2, y_3\right\},\left\{y_1, y_2\right\} \rightarrow x_3 \rightarrow\left\{y_3, y_4\right\}$.

We know $y_1 \in V_2^{+--}, y_2 \in V_2^{-+-}, y_3 \in V_2^{-++}$and $y_4 \in V_2^{--+}$.  \ By Lemma 2.1, we can get $V_3^{+--}=V_3^{-+-}=V_3^{-++}=V_3^{--+}=\varnothing$. \ Thus $V_3=V_3^{+++} \cup V_3^{++-} \cup V_3^{+-+} \cup V_3^{---}$. \ Since $\partial_D\left(V_3^{++-}, y_i\right) \leq 2$ where $i=1,2,3$, then we have $V_3^{++-} \rightarrow\left\{y_1, y_2\right\}$. \ Since $\partial_D\left(V_3^{+-+}, y_i\right) \leq 2$ where $i=1,2,3$, then we have $V_3^{+-+} \rightarrow\left\{y_1, y_4\right\}$. \ Let $\left|V_3^{++-}\right|=q_1,\left|V_3^{+-+}\right|=q_2$ and $F_1=D\left[V_2 \backslash\left\{y_1, y_2\right\} \cup V_3^{++-}\right], F_2=D\left[V_2 \backslash\left\{y_1, y_4\right\} \cup V_3^{+-+}\right]$. \ Then $F_1$ and $F_2$ are respextively an orientation of $K\left(2, q_1\right)$ and $K\left(2, q_2\right)$ where $q_1, q_2 \leq q$. \ If $q_1 \geq 3$, then there exists $z_i, z_j \in V_3^{++-}$such that $\partial_{F_1}\left(z_i, z_j\right)=4$, so $\partial_D\left(z_i, z_j\right) \geq 3$, a contradiction. \ Hence $q_1 \leq 2$. \ If $q_2 \geq 3$, then there exists $z_i, z_j \in V_3^{+-+}$such that $\partial_{F_2}\left(z_i, z_j\right)=4$, so $\partial_D\left(z_i, z_j\right) \geq 3$, a contradiction. \ Hence $q_2 \leq 2$. \ By Lemma 2.1, we also have $\left|V_3^{+++}\right| \leq 1$ and $\left|V_3^{---}\right| \leq 1$, thus $\left|V_3\right|=q \leq 1+2+2+1=6 \leq 11$.\\
(18) Case $(1,2,3)$.\\
Subcase 1: $V_2 \backslash\left\{y_1\right\} \rightarrow x_1 \rightarrow y_1, V_2 \backslash\left\{y_1, y_2\right\} \rightarrow x_2 \rightarrow\left\{y_1, y_2\right\}, y_4 \rightarrow x_3 \rightarrow V_2 \backslash\left\{y_4\right\}$.

We know $y_1 \in V_2^{+++}, y_2 \in V_2^{-++}, y_3 \in V_2^{--+}$and $y_4 \in V_2^{---}$. By Lemma 2.1, we can get $V_3^{+++}=V_3^{-++}=V_3^{--+}=V_3^{---}=\varnothing$. \ Thus $V_3=V_3^{++-} \cup V_3^{+-+} \cup V_3^{+--} \cup V_3^{-+-}$. \ Take any $z \in V_3$, since $\partial_D\left(y_1, z\right) \leq 2$ and $\partial_D\left(z, y_4\right) \leq 2$, then we have $y_1 \rightarrow z$ and $z \rightarrow y_4$, this means $y_1 \rightarrow V_3$ and $V_3 \rightarrow y_4$. \ Since $\partial_D\left(V_3^{+-+}, y_i\right) \leq 2$ where $i=1,2,3$, then we have $V_3^{+-+} \rightarrow y_3$. \ Since $\partial_D\left(y_i, V_3^{-+-}\right) \leq 2$ where $i=1,2,3$, then we have $y_2 \rightarrow V_3^{-+-}$.  \ Let $\left|V_3^{++-}\right|=q_1,\left|V_3^{+-+}\right|=q_2,\left|V_3^{+--}\right|=q_3,\left|V_3^{-+-}\right|=q_4$ and $F_1=D\left[V_2 \backslash\left\{y_1, y_4\right\} \cup V_3^{++-}\right], F_2=D\left[V_2 \backslash\left\{y_1, y_4\right\} \cup V_3^{+--}\right]$.  \ Then $F_1$ and $F_2$ are respextively an orientation of $K\left(2, q_1\right)$ and $K\left(2, q_3\right)$ where $q_1, q_3 \leq q$. \ If $q_1 \geq 3$, then there exists $z_i, z_j \in V_3^{++-}$such that $\partial_{F_1}\left(z_i, z_j\right)=4$, so $\partial_D\left(z_i, z_j\right) \geq 3$, a contradiction. \ Hence $q_1 \leq 2$. \ If $q_2 \geq 2$, then there exists $z_i, z_j \in V_3^{+-+}$such that $\partial_D\left(z_i, z_j\right) \geq 3$, a contradiction. \ Hence $q_2 \leq 1$. \ If $q_3 \geq 3$, then there exists $z_i, z_j \in V_3^{+--}$such that $\partial_{F_2}\left(z_i, z_j\right)=4$, so $\partial_D\left(z_i, z_j\right) \geq 3$, a contradiction. \ Hence $q_3 \leq 2$. \ If $q_4 \geq 2$, then there exists $z_i, z_j \in V_3^{-+-}$such that $\partial_D\left(z_i, z_j\right) \geq 3$, a contradiction. \ Hence $q_4 \leq 1$. \ Thus $\left|V_3\right|=q \leq 2+1+2+1=6 \leq 11$.\\
Subcase 2: $V_2 \backslash\left\{y_1\right\} \rightarrow x_1 \rightarrow y_1, V_2 \backslash\left\{y_1, y_2\right\} \rightarrow x_2 \rightarrow\left\{y_1, y_2\right\}, y_2 \rightarrow x_3 \rightarrow V_2 \backslash\left\{y_2\right\}$.

We know $y_1 \in V_2^{+++}, y_2 \in V_2^{-+-}$and $y_3, y_4 \in V_2^{--+}$. \ By Lemma 2.1, we can get $V_3^{+++}=V_3^{-+-}=V_3^{--+}=\varnothing$. \ Thus $V_3=V_3^{++-} \cup V_3^{+-+} \cup V_3^{-++} \cup V_3^{+--} \cup V_3^{---}$. \ Take any $z \in V_3$, since $\partial_D\left(y_1, z\right) \leq 2$, then we have $y_1 \rightarrow z$, this means $y_1 \rightarrow V_3$. \ Since $\partial_D\left(V_3^{++-}, y_i\right) \leq 2$ where $i=1,2,3$, then we have $V_3^{++-} \rightarrow y_2$. \ Since $\partial_D\left(V_3^{+-+}, y_i\right) \leq 2$ where $i=1,2,3$, then we have $V_3^{+-+} \rightarrow\left\{y_3, y_4\right\}$.  \ Since $\partial_D\left(V_3^{-++}, y_i\right) \leq 2$ where $i=1,2,3$, then we have $V_3^{-++} \rightarrow\left\{y_2, y_3, y_4\right\}$. \ Let $\left|V_3^{++-}\right|=q_1,\left|V_3^{+-+}\right|=q_2,\left|V_3^{-++}\right|=q_3,\left|V_3^{+--}\right|=q_4$ and $F_1=D\left[V_2 \backslash\left\{y_1, y_2\right\} \cup V_3^{++-}\right], F_2=D\left[V_2 \backslash\left\{y_1\right\} \cup V_3^{+--}\right]$. \ Then $F_1$ and $F_2$ are respextively an orientation of $K\left(2, q_1\right)$ and $K\left(3, q_4\right)$ where $q_1, q_4 \leq q$. \ If $q_1 \geq 3$, then there exists $z_i, z_j \in V_3^{++-}$such that $\partial_{F_1}\left(z_i, z_j\right)=4$, so $\partial_D\left(z_i, z_j\right) \geq 3$, a contradiction. \ Hence $q_1 \leq 2$. \ If $q_2 \geq 2$, then there exists $z_i, z_j \in V_3^{+-+}$such that $\partial_D\left(z_i, z_j\right) \geq 3$, a contradiction. \ Hence $q_2 \leq 1$. \ If $q_3 \geq 2$, then there exists $z_i, z_j \in V_3^{-++}$such that $\partial_D\left(z_i, z_j\right) \geq 3$, a contradiction. \ Hence $q_3 \leq 1$. \ If $q_4 \geq 4$, then there exists $z_i, z_j \in V_3^{+--}$such that $\partial_{F_2}\left(z_i, z_j\right)=4$, so $\partial_D\left(z_i, z_j\right) \geq 3$, a contradiction. \ Hence $q_4 \leq 3$. \ By Lemma 2.1, we also have $\left|V_3^{---}\right| \leq 1$, \ thus $\left|V_3\right|=q \leq 2+1+1+3+1=8 \leq 11$.\\
Subcase 3: $V_2 \backslash\left\{y_1\right\} \rightarrow x_1 \rightarrow y_1, V_2 \backslash\left\{y_1, y_2\right\} \rightarrow x_2 \rightarrow\left\{y_1, y_2\right\}, y_1 \rightarrow x_3 \rightarrow V_2 \backslash\left\{y_1\right\}$.

We know $y_1 \in V_2^{++-}, y_2 \in V_2^{-++}$and $y_3, y_4 \in V_2^{--+}$. \ By Lemma 2.1, we can get $V_3^{++-}=V_3^{-++}=V_3^{--+}=\varnothing$. \ Thus $V_3=V_3^{+++} \cup V_3^{+-+} \cup V_3^{+--} \cup V_3^{-+-} \cup V_3^{---}$. \ Since $\partial_D\left(V_3^{+-+}, y_i\right) \leq 2$ where $i=1,2,3$, then we have $V_3^{+-+} \rightarrow\left\{y_3, y_4\right\}$. \ Since $\partial_D\left(y_i, V_3^{+--}\right) \leq 2$ where $i=1,2,3$, then we have $y_1 \rightarrow V_3^{+--}$. \ Since $\partial_D\left(y_i, V_3^{-+-}\right) \leq 2$ where $i=1,2,3$, then we have $\left\{y_1, y_2\right\} \rightarrow V_3^{-+-}$. \ Let $\left|V_3^{+-+}\right|=q_1,\left|V_3^{+--}\right|=q_2,\left|V_3^{-+-}\right|=q_3$ and $
F_1=D\left[V_2 \backslash\left\{y_3, y_4\right\} \cup V_3^{+-+}\right], F_2=D\left[V_2 \backslash\left\{y_1\right\} \cup V_3^{+--}\right], F_3=D\left[V_2 \backslash\left\{y_1, y_2\right\} \cup V_3^{-+-}\right] .$  \ Then $F_1, F_2$ and $F_3$ are respextively an orientation of $K\left(2, q_1\right), K\left(3, q_2\right)$ and $K\left(2, q_3\right)$ where $q_1, q_2, q_3 \leq q$.  \ If $q_1 \geq 3$, then there exists $z_i, z_j \in V_3^{+-+}$such that $\partial_{F_1}\left(z_i, z_j\right)=4$, so $\partial_D\left(z_i, z_j\right) \geq 3$, a contradiction. \ Hence $q_1 \leq 2$. \ If $q_2 \geq 4$, then there exists $z_i, z_j \in V_3^{+--}$such that $\partial_{F_2}\left(z_i, z_j\right)=4$, so $\partial_D\left(z_i, z_j\right) \geq 3$, a contradiction. \ Hence $q_2 \leq 3$. \ If $q_3 \geq 3$, then there exists $z_i, z_j \in V_3^{-+-}$such that $\partial_{F_3}\left(z_i, z_j\right)=4$, so $\partial_D\left(z_i, z_j\right) \geq 3$, a contradiction. \ Hence $q_3 \leq 2$. \ By Lemma 2.1, we also have $\left|V_3^{+++}\right| \leq 1$ and $\left|V_3^{---}\right| \leq 1$, \ thus $\left|V_3\right|=q \leq 1+2+3+2+1=9 \leq 11$.\\
Subcase 4: $V_2 \backslash\left\{y_1\right\} \rightarrow x_1 \rightarrow y_1, V_2 \backslash\left\{y_2, y_3\right\} \rightarrow x_2 \rightarrow\left\{y_2, y_3\right\}, y_4 \rightarrow x_3 \rightarrow V_2 \backslash\left\{y_4\right\}$.

We know $y_1 \in V_2^{+-+}, y_2, y_3 \in V_2^{-++}$and $y_4 \in V_2^{---}$. \ By Lemma 2.1, we can get $V_3^{+-+}=V_3^{-++}=V_3^{---}=\varnothing$. \ Thus $V_3=V_3^{+++} \cup V_3^{++-} \cup V_3^{+--} \cup V_3^{-+-} \cup V_3^{--+}$. \ Take any $z \in V_3$, since $\partial_D\left(z, y_4\right) \leq 2$, then we have $z \rightarrow y_4$, this means $V_3 \rightarrow y_4$. \ Since $\partial_D\left(y_i, V_3^{+--}\right) \leq 2$ where $i=1,2,3$, then we have $y_1 \rightarrow V_3^{+--}$. \ Since $\partial_D\left(y_i, V_3^{-+-}\right) \leq 2$ where $i=1,2,3$, then we have $\left\{y_2, y_3\right\} \rightarrow V_3^{-+-}$.  \ Since $\partial_D\left(y_i, V_3^{--+}\right) \leq 2$ where $i=1,2,3$, then we have $\left\{y_1, y_2, y_3\right\} \rightarrow V_3^{--+}$. \ Let $\left|V_3^{++-}\right|=q_1,\left|V_3^{+--}\right|=q_2,\left|V_3^{-+-}\right|=q_3,\left|V_3^{--+}\right|=q_4$ and $F_1=D\left[V_2 \backslash\left\{y_4\right\} \cup V_3^{++-}\right], F_2=D\left[V_2 \backslash\left\{y_1, y_4\right\} \cup V_3^{+--}\right]$. \ Then $F_1$ and $F_2$ are respextively an orientation of $K\left(3, q_1\right)$ and $K\left(2, q_2\right)$ where $q_1, q_2 \leq q$. \ If $q_1 \geq 4$, then there exists $z_i, z_j \in V_3^{++-}$such that $\partial_{F_1}\left(z_i, z_j\right)=4$, so $\partial_D\left(z_i, z_j\right) \geq 3$, a contradiction. \ Hence $q_1 \leq 3$. \ If $q_2 \geq 3$, then there exists $z_i, z_j \in V_3^{+--}$such that $\partial_{F_2}\left(z_i, z_j\right)=4$, so $\partial_D\left(z_i, z_j\right) \geq 3$, a contradiction. \ Hence $q_2 \leq 2$. \ If $q_3 \geq 2$, then there exists $z_i, z_j \in V_3^{-+-}$such that $\partial_D\left(z_i, z_j\right) \geq 3$, a contradiction. \ Hence $q_3 \leq 1$.\ If $q_4 \geq 2$, then there exists $z_i, z_j \in V_3^{--+}$such that $\partial_D\left(z_i, z_j\right) \geq 3$, a contradiction. \ Hence $q_4 \leq 1$. \ By Lemma 2.1, we also have $\left|V_3^{+++}\right| \leq 1$, \ thus $\left|V_3\right|=q \leq 1+3+2+1+1=8 \leq 11$.\\
Subcase 5: $V_2 \backslash\left\{y_1\right\} \rightarrow x_1 \rightarrow y_1, V_2 \backslash\left\{y_2, y_3\right\} \rightarrow x_2 \rightarrow\left\{y_2, y_3\right\}, y_2 \rightarrow x_3 \rightarrow V_2 \backslash\left\{y_2\right\}$.

We know $y_1 \in V_2^{+-+}, y_2 \in V_2^{-+-}, y_3 \in V_2^{-++}$and $y_4 \in V_2^{--+}$. \ By Lemma 2.1, we can get $V_3^{+-+}=V_3^{-+-}=V_3^{-++}=V_3^{--+}=\varnothing$. \ Thus $V_3=V_3^{+++} \cup V_3^{++-} \cup V_3^{+--} \cup V_3^{---}$. \ Since $\partial_D\left(V_3^{++-}, y_i\right) \leq 2$ where $i=1,2,3$, then we have $V_3^{++-} \rightarrow y_2$. \ Since $\partial_D\left(y_i, V_3^{+--}\right) \leq 2$ where $i=1,2,3$, then we have $y_1 \rightarrow V_3^{+--}$. \ Let $\left|V_3^{++-}\right|=q_1,\left|V_3^{+--}\right|=q_2$ and $F_1=D\left[V_2 \backslash\left\{y_2\right\} \cup V_3^{++-}\right], F_2=D\left[V_2 \backslash\left\{y_1\right\} \cup V_3^{+--}\right]$. \ Then $F_1$ and $F_2$ are respextively an orientation of $K\left(3, q_1\right)$ and $K\left(3, q_2\right)$ where $q_1, q_2 \leq q$. \ If $q_1 \geq 4$, then there exists $z_i, z_j \in V_3^{++-}$such that $\partial_{F_1}\left(z_i, z_j\right)=4$, so $\partial_D\left(z_i, z_j\right) \geq 3$, a contradiction. \ Hence $q_1 \leq 3$. \ If $q_2 \geq 4$, then there exists $z_i, z_j \in V_3^{+--}$such that $\partial_{F_2}\left(z_i, z_j\right)=4$, so $\partial_D\left(z_i, z_j\right) \geq 3$, a contradiction. \ Hence $q_2 \leq 3$. \ By Lemma 2.1, we also have $\left|V_3^{+++}\right| \leq 1$ and $\left|V_3^{---}\right| \leq 1$, thus $\left|V_3\right|=q \leq 1+3+3+1=8 \leq 11$.\\
Subcase 6: $V_2 \backslash\left\{y_1\right\} \rightarrow x_1 \rightarrow y_1, V_2 \backslash\left\{y_2, y_3\right\} \rightarrow x_2 \rightarrow\left\{y_2, y_3\right\}, y_1 \rightarrow x_3 \rightarrow V_2 \backslash\left\{y_1\right\}$.

We know $y_1 \in V_2^{+--}, y_2, y_3 \in V_2^{-++}$and $y_4 \in V_2^{--+}$. \ By Lemma 2.1, we can get $V_3^{+--}=V_3^{-++}=V_3^{--+}=\varnothing$. \ Thus $V_3=V_3^{+++} \cup V_3^{++-} \cup V_3^{+-+} \cup V_3^{-+-} \cup V_3^{---}$. \ Since $\partial_D\left(V_3^{++-}, y_i\right) \leq 2$ where $i=1,2,3$, then we have $V_3^{++-} \rightarrow y_1$. \ Since $\partial_D\left(V_3^{+-+}, y_i\right) \leq 2$ where $i=1,2,3$, then we have $V_3^{+-+} \rightarrow\left\{y_1, y_4\right\}$. \ Since $\partial_D\left(y_i, V_3^{-+-}\right) \leq 2$ where $i=1,2,3$, then we have $\left\{y_2, y_3\right\} \rightarrow V_3^{-+-}$. \ Let $\left|V_3^{++-}\right|=q_1,\left|V_3^{+-+}\right|=q_2,\left|V_3^{-+-}\right|=q_3$ and $
F_1=D\left[V_2 \backslash\left\{y_1\right\} \cup V_3^{++-}\right], F_2=D\left[V_2 \backslash\left\{y_1, y_4\right\} \cup V_3^{+-+}\right], F_3=D\left[V_2 \backslash\left\{y_2, y_3\right\} \cup V_3^{-+-}\right] .$ \ Then $F_1, F_2$ and $F_3$ are respextively an orientation of $K\left(3, q_1\right), K\left(2, q_2\right)$ and $K\left(2, q_3\right)$ where $q_1, q_2, q_3 \leq q$. \ If $q_1 \geq 4$, then there exists $z_i, z_j \in V_3^{++-}$such that $\partial_{F_1}\left(z_i, z_j\right)=4$, so $\partial_D\left(z_i, z_j\right) \geq 3$, a contradiction. \ Hence $q_1 \leq 3$. \ If $q_2 \geq 3$, then there exists $z_i, z_j \in V_3^{+-+}$such that $\partial_{F_2}\left(z_i, z_j\right)=4$, so $\partial_D\left(z_i, z_j\right) \geq 3$, a contradiction. \ Hence $q_2 \leq 2$. \ If $q_3 \geq 3$, then there exists $z_i, z_j \in V_3^{-+-}$such that $\partial_{F_3}\left(z_i, z_j\right)=4$, so $\partial_D\left(z_i, z_j\right) \geq 3$, a contradiction. \ Hence $q_3 \leq 2$. \ By Lemma 2.1, we also have $\left|V_3^{+++}\right| \leq 1$ and $\left|V_3^{---}\right| \leq 1$, \ thus $\left|V_3\right|=q \leq 1+3+2+2+1=9 \leq 11$.\\
(19) Case $(2,2,2)$.\\
Subcase 1: $V_2 \backslash\left\{y_1, y_2\right\} \rightarrow V_1 \rightarrow\left\{y_1, y_2\right\}$.

Take any $y \in\left\{y_1, y_2\right\}$ and $z \in V_3$, since $\partial_D(y, z) \leq 2$, then we have $y \rightarrow z$, this means $\left\{y_1, y_2\right\} \rightarrow V_3$. So $\partial_D\left(y_1, y_2\right) \geq 3$, a contradiction.\\
Subcase 2: $\left\{y_3, y_4\right\} \rightarrow\left\{x_1, x_2\right\} \rightarrow\left\{y_1, y_2\right\},\left\{y_2, y_4\right\} \rightarrow x_3 \rightarrow\left\{y_1, y_3\right\}$.

We know $y_1 \in V_2^{+++}, y_2 \in V_2^{++-}, y_3 \in V_2^{--+}$and $y_4 \in V_2^{---}$.  \ By Lemma 2.1, we can get $V_3^{+++}=V_3^{++-}=V_3^{--+}=V_3^{---}=\varnothing$. \ Thus $V_3=V_3^{+-+} \cup V_3^{-++} \cup V_3^{+--} \cup V_3^{-+-}$. \ Take any $z \in V_3$, since $\partial_D\left(y_1, z\right) \leq 2$ and $\partial_D\left(z, y_4\right) \leq 2$, then we have $y_1 \rightarrow z$ and $z \rightarrow y_4$, this means $y_1 \rightarrow V_3$ and $V_3 \rightarrow y_4$. \ Since $\partial_D\left(V_3^{+-+}, y_i\right) \leq 2$ where $i=1,2,3$, then we have $V_3^{+-+} \rightarrow y_3$. \ Since $\partial_D\left(V_3^{-++}, y_i\right) \leq 2$ where $i=1,2,3$, then we have $V_3^{-++} \rightarrow y_3$. \ Since $\partial_D\left(y_i, V_3^{+--}\right) \leq 2$ where $i=1,2,3$, then we have $y_2 \rightarrow V_3^{+--}$. \ Since $\partial_D\left(y_i, V_3^{-+-}\right) \leq 2$ where $i=1,2,3$, then we have $y_2 \rightarrow V_3^{-+-}$. \ Let $\left|V_3^{+-+}\right|=q_1,\left|V_3^{-++}\right|=q_2,\left|V_3^{+--}\right|=q_3$ and $\left|V_3^{-+-}\right|=q_4$. \ If $q_1 \geq 2$, then there exists $z_i, z_j \in V_3^{+-+}$such that $\partial_D\left(z_i, z_j\right) \geq 3$, a contradiction.
\ Hence $q_1 \leq 1$. \ If $q_2 \geq 2$, then there exists $z_i, z_j \in V_3^{-++}$such that $\partial_D\left(z_i, z_j\right) \geq 3$, a contradiction. \ Hence $q_2 \leq 1$. \ If $q_3 \geq 2$, then there exists $z_i, z_j \in V_3^{+--}$such that $\partial_D\left(z_i, z_j\right) \geq 3$, a contradiction. \ Hence $q_3 \leq 1$. \ If $q_4 \geq 2$, then there exists $z_i, z_j \in V_3^{-+-}$such that $\partial_D\left(z_i, z_j\right) \geq 3$, a contradiction. \ Hence $q_4 \leq 1$. \ Thus $\left|V_3\right|=q \leq 1+1+1+1=4 \leq 11$.\\
Subcase 3: $\left\{y_3, y_4\right\} \rightarrow\left\{x_1, x_2\right\} \rightarrow\left\{y_1, y_2\right\} \rightarrow x_3 \rightarrow\left\{y_3, y_4\right\}$.

We know $y_1, y_2 \in V_2^{++-}$and $y_3, y_4 \in V_2^{--+}$. \ By Lemma 2.1, we can get $V_3^{++-}=V_3^{--+}=\varnothing$. \ Thus $V_3=V_3^{+++} \cup V_3^{+-+} \cup V_3^{-++} \cup V_3^{+--} \cup V_3^{-+-} \cup V_3^{---}$. \ Since $\partial_D\left(V_3^{+-+}, y_i\right) \leq 2$ where $i=1,2,3$, then we have $V_3^{+-+} \rightarrow\left\{y_3, y_4\right\}$. \ Since $\partial_D\left(V_3^{-++}, y_i\right) \leq 2$ where $i=1,2,3$, then we have $V_3^{-++} \rightarrow\left\{y_3, y_4\right\}$. \ Since $\partial_D\left(y_i, V_3^{+--}\right) \leq 2$ where $i=1,2,3$, then we have $\left\{y_1, y_2\right\} \rightarrow V_3^{+--}$. \ Since $\partial_D\left(y_i, V_3^{-+-}\right) \leq 2$ where $i=1,2,3$, then we have $\left\{y_1, y_2\right\} \rightarrow V_3^{-+-}$. \ Let $\left|V_3^{+-+}\right|=q_1,\left|V_3^{-++}\right|=q_2,\left|V_3^{+--}\right|=q_3,\left|V_3^{-+-}\right|=q_4$ and $F_1=D\left[V_2^{++-} \cup V_3^{+-+}\right], F_2=D\left[V_2^{++-} \cup V_3^{-++}\right], \allowbreak F_3=D\left[V_2^{--+} \cup V_3^{+--}\right], F_4=D\left[V_2^{--+} \cup V_3^{-+-}\right]$.
\ Then $F_1, F_2, F_3$ and $F_4$ are respextively an orientation of $K\left(2, q_1\right), K\left(2, q_2\right), K\left(2, q_3\right)$ and $K\left(2, q_4\right)$ where $q_1, q_2, q_3, q_4 \leq q$. \ If $q_1 \geq 3$, then there exists $z_i, z_j \in V_3^{+-+}$such that $\partial_{F_1}\left(z_i, z_j\right)=4$, so $\partial_D\left(z_i, z_j\right) \geq 3$, a contradiction. \ Hence $q_1 \leq 2$. \ If $q_2 \geq 3$, then there exists $z_i, z_j \in V_3^{-++}$such that $\partial_{F_2}\left(z_i, z_j\right)=4$, so $\partial_D\left(z_i, z_j\right) \geq 3$, a contradiction. \ Hence $q_2 \leq 2$. \ If $q_3 \geq 3$, then there exists $z_i, z_j \in V_3^{+--}$such that $\partial_{F_3}\left(z_i, z_j\right)=4$, so $\partial_D\left(z_i, z_j\right) \geq 3$, a contradiction. \ Hence $q_3 \leq 2$. \ If $q_4 \geq 3$, then there exists $z_i, z_j \in V_3^{-+-}$such that $\partial_{F_4}\left(z_i, z_j\right)=4$, so $\partial_D\left(z_i, z_j\right) \geq 3$, a contradiction. \ Hence $q_4 \leq 2$. \ By Lemma 2.1, we also have $\left|V_3^{+++}\right| \leq 1$ and $\left|V_3^{---}\right| \leq 1$, \ thus $\left|V_3\right|=q \leq 1+2+2+2+2+1=10 \leq 11$.\\
Subcase 4: $\left\{y_3, y_4\right\} \rightarrow x_1 \rightarrow\left\{y_1, y_2\right\},\left\{y_2, y_4\right\} \rightarrow x_2 \rightarrow\left\{y_1, y_3\right\},\left\{y_2, y_3\right\} \rightarrow x_3 \rightarrow\left\{y_1, y_4\right\}$.

We know $y_1 \in V_2^{+++}, y_2 \in V_2^{+--}, y_3 \in V_2^{-+-}$and $y_4 \in V_2^{--+}$.  \ By Lemma 2.1, we can get $V_3^{+++}=V_3^{+--}=V_3^{-+-}=V_3^{--+}=\varnothing$. \ Thus $V_3=V_3^{++-} \cup V_3^{+-+} \cup V_3^{-++} \cup V_3^{---}$. \ Take any $z \in V_3$, since $\partial_D\left(y_1, z\right) \leq 2$, then we have $y_1 \rightarrow z$, this means $y_1 \rightarrow V_3$. \ Since $\partial_D\left(V_3^{++-}, y_i\right) \leq 2$ where $i=1,2,3$, then we have $V_3^{++-} \rightarrow\left\{y_2, y_3\right\}$. \ Since $\partial_D\left(V_3^{+-+}, y_i\right) \leq 2$ where $i=1,2,3$, then we have $V_3^{+-+} \rightarrow\left\{y_2, y_4\right\}$. \ Since $\partial_D\left(V_3^{-++}, y_i\right) \leq 2$ where $i=1,2,3$, then we have $V_3^{-++} \rightarrow\left\{y_3, y_4\right\}$. \ Let $\left|V_3^{++-}\right|=q_1,\left|V_3^{+-+}\right|=q_2$ and $\left|V_3^{-++}\right|=q_3$.  \ If $q_1 \geq 2$, then there exists $z_i, z_j \in V_3^{++-}$ such that $\partial_D\left(z_i, z_j\right) \geq 3$, a contradiction. \ Hence $q_1 \leq 1$. \ If $q_2 \geq 2$, then there exists $z_i, z_j \in V_3^{+-+}$such that $\partial_D\left(z_i, z_j\right) \geq 3$, a contradiction. \ Hence $q_2 \leq 1$. \ If $q_3 \geq 2$, then there exists $z_i, z_j \in V_3^{-++}$such that $\partial_D\left(z_i, z_j\right) \geq 3$, a contradiction. \ Hence $q_3 \leq 1$. \ By Lemma 2.1, we also have $\left|V_3^{---}\right| \leq 1$, \ thus $\left|V_3\right|=q \leq 1+1+1+1=4 \leq 11$.\\
Subcase 5: $\left\{y_3, y_4\right\} \rightarrow x_1 \rightarrow\left\{y_1, y_2\right\},\left\{y_2, y_4\right\} \rightarrow x_2 \rightarrow\left\{y_1, y_3\right\},\left\{y_1, y_4\right\} \rightarrow x_3 \rightarrow\left\{y_2, y_3\right\}$.

We know $y_1 \in V_2^{++-}, y_2 \in V_2^{+-+}, y_3 \in V_2^{-++}$and $y_4 \in V_2^{---}$.  \ By Lemma 2.1, we can get $V_3^{++-}=V_3^{+-+}=V_3^{-++}=V_3^{---}=\varnothing$. \ Thus $V_3=V_3^{+++} \cup V_3^{+--} \cup V_3^{-+-} \cup V_3^{--+}$. \ Take any $z \in V_3$, since $\partial_D\left(z, y_4\right) \leq 2$, then we have $z \rightarrow y_4$, this means $V_3 \rightarrow y_4$. \ Since $\partial_D\left(y_i, V_3^{+--}\right) \leq 2$ where $i=1,2,3$, then we have $\left\{y_1, y_2\right\} \rightarrow V_3^{+--}$. \ Since $\partial_D\left(y_i, V_3^{-+-}\right) \leq 2$ where $i=1,2,3$, then we have $\left\{y_1, y_3\right\} \rightarrow V_3^{-+-}$. \ Since $\partial_D\left(y_i, V_3^{--+}\right) \leq 2$ where $i=1,2,3$, then we have $\left\{y_2, y_3\right\} \rightarrow V_3^{--+}$. \ Let $\left|V_3^{+--}\right|=q_1,\left|V_3^{-+-}\right|=q_2$ and $\left|V_3^{--+}\right|=q_3$.  \ If $q_1 \geq 2$, then there exists $z_i, z_j \in V_3^{+--}$ such that $\partial_D\left(z_i, z_j\right) \geq 3$, a contradiction. \ Hence $q_1 \leq 1$. \ If $q_2 \geq 2$, then there exists $z_i, z_j \in V_3^{-+-}$such that $\partial_D\left(z_i, z_j\right) \geq 3$, a contradiction. \ Hence $q_2 \leq 1$. \ If $q_3 \geq 2$, then there exists $z_i, z_j \in V_3^{--+}$such that $\partial_D\left(z_i, z_j\right) \geq 3$, a contradiction. \ Hence $q_3 \leq 1$. \ By Lemma 2.1, we also have $\left|V_3^{+++}\right| \leq 1$, \ thus $\left|V_3\right|=q \leq 1+1+1+1=4 \leq 11$.\\
Subcase $6:\left\{y_3, y_4\right\} \rightarrow x_1 \rightarrow\left\{y_1, y_2\right\},\left\{y_2, y_4\right\} \rightarrow x_2 \rightarrow\left\{y_1, y_3\right\},\left\{y_1, y_2\right\} \rightarrow x_3 \rightarrow\left\{y_3, y_4\right\}$.

We know $y_1 \in V_2^{++-}, y_2 \in V_2^{+--}, y_3 \in V_2^{-++}$and $y_4 \in V_2^{--+}$.  \ By Lemma 2.1, we can get $V_3^{++-}=V_3^{+--}=V_3^{-++}=V_3^{--+}=\varnothing$. \ Thus $V_3=V_3^{+++} \cup V_3^{+-+} \cup V_3^{-+-} \cup V_3^{---}$. \ Since $\partial_D\left(V_3^{+-+}, y_i\right) \leq 2$ where $i=1,2,3$, then we have $V_3^{+-+} \rightarrow\left\{y_2, y_4\right\}$.  \ Since $\partial_D\left(y_i, V_3^{-+-}\right) \leq 2$ where $i=1,2,3$, then we have $\left\{y_1, y_3\right\} \rightarrow V_3^{-+-}$. \ Let $\left|V_3^{+-+}\right|=q_1,\left|V_3^{-+-}\right|=q_2$ and
$F_1=D\left[V_2 \backslash\left\{y_2, y_4\right\} \cup V_3^{+-+}\right], F_2=D\left[V_2 \backslash\left\{y_1, y_3\right\} \cup V_3^{-+-}\right]$. \ Then $F_1$ and $F_2$ are respextively an orientation of $K\left(2, q_1\right)$ and $K\left(2, q_2\right)$ where $q_1, q_2 \leq q$. \ If $q_1 \geq 3$, then there exists $z_i, z_j \in V_3^{+-+}$such that $\partial_{F_1}\left(z_i, z_j\right)=4$, so $\partial_D\left(z_i, z_j\right) \geq 3$, a contradiction. \ Hence $q_1 \leq 2$. \ If $q_2 \geq 3$, then there exists $z_i, z_j \in V_3^{-+-}$such that $\partial_{F_2}\left(z_i, z_j\right)=4$, so $\partial_D\left(z_i, z_j\right) \geq 3$, a contradiction. \ Hence $q_2 \leq 2$. \ By Lemma 2.1, we also have $\left|V_3^{+++}\right| \leq 1$ and $\left|V_3^{---}\right| \leq 1$, \ thus $\left|V_3\right|=q \leq 1+2+2+1=6 \leq 11$.\\

In summary,\ it can be concluded that if $f(K(3,4, q))=2$, then $q \leq 11$. \  Since when $q \leq 11$, we have found an orientation of diameter 2 of $K(3,4, q)$ in Lemma 4.1. \ Therefore, $f(K(3,4, q))=2$ if and only if $q \leq 11$.   $\hfill\blacksquare$\\

\bibliography{sn-bibliography}% common bib file
%\bibliographystyle{unsrt}
%% if required, the content of .bbl file can be included here once bbl is generated
%%\input sn-article.bbl

\end{document}